\documentclass[reqno,a4paper]{amsart}

\usepackage{amsmath}
\usepackage{amssymb}
\usepackage{amsthm}
\usepackage{enumerate}
\usepackage{mathrsfs} 
\usepackage{hyperref}

\usepackage{eqlist}
\usepackage{array}

\theoremstyle{plain}
\newtheorem{theorem}{Theorem}[section]

\newtheorem{prop}[theorem]{Proposition}

\newtheorem{lemma}{Lemma}[section]
\newtheorem{corollary}{Corollary}[theorem]

\theoremstyle{definition}

\newtheorem{example}{\textit{Example}}

\numberwithin{equation}{section}

\newcommand{\be}{\begin{equation}}
\newcommand{\ee}{\end{equation}}

\begin{document}

\title[Parabolic frequency on compact Finsler metric measure spaces]
{Parabolic frequency and backward uniqueness on compact Finsler metric measure spaces}
\author[Bin Shen \and Zisu Zhao]
{Bin Shen* \and Zisu Zhao**}

\newcommand{\acr}{\newline\indent}

\address{\llap{*\,}School of Mathematics\acr
                   Southeast University\acr
                   Dongnandaxue Road 2\acr
                   Nanjing, 211189, Jiangsu\acr
                   China}
\email{shenbin@seu.edu.cn}

\address{\llap{**\,}Institute of Mathematics\acr
		Academy of Mathematics and Systems Science\acr
		Chinese Academy of Sciences\acr
		Zhongguancun East Road 55\acr
		100190, Beijing\acr
		China}
\email{steven1998zhao@outlook.com}

\thanks{This work was partially supported by the National Natural Science Foundation of China (Grant Nos. 12001099 and 12271093).}


\begin{abstract}
We establish monotonicity formulas for parabolic frequency functions on closed Finsler metric measure spaces with fixed or evolving metrics and measures. When the reference measure is fixed, an integral bound on the metric variation yields monotonicity for the Finsler heat equation. For coupled metric measure flows, we construct positive global supersolutions and subsolutions of a linearized conjugate heat equation from local radial barriers. These weights define a frequency function whose monotonicity follows from an integral bound involving the metric variation and the change in weighted Ricci curvature. For parabolic equations with a source controlled by the solution and its gradient, we obtain a differential inequality for the frequency and prove that a solution vanishing at the terminal time must vanish throughout the preceding interval. The proof controls the comparison of the two weights through an energy quotient computed with the subsolution measure.
\end{abstract}

\maketitle
\enlargethispage{4pt}

\tableofcontents

\section{Introduction}

In recent years, geometric analysis on Finsler metric measure spaces has developed rapidly. In the theory of heat equations, S. Ohta and K.-T. Sturm characterized Finslerian heat flows as $L^{2}$-gradient flows, either of the energy functional on the manifold or of the relative entropy on the reverse $L^{2}$-Wasserstein space. They proved the equivalence of these two formulations and obtained regularity results for solutions, together with pointwise comparison estimates and integrated upper Gaussian bounds \cite{OhtaSturm2009HeatFlow}. Subsequently, Q. Xia proved that nonnegative $L^{p}$-global subsolutions to the heat equation on forward complete Finsler metric measure spaces are uniquely determined by their initial data \cite{Xia2024LpSubsolutions}. These works show that Finslerian parabolic equations admit a well-developed theory of forward well-posedness and uniqueness. Nevertheless, compared with the corresponding theory on Riemannian manifolds \cite{Ni2015ParabolicFrequency, LiWang2018Parabolic}, several fundamental problems for Finslerian heat equations remain insufficiently understood, especially backward uniqueness, quantitative unique continuation, and vanishing-order estimates.

The parabolic frequency function is a central tool in the study of such problems in Riemannian geometric analysis. The method is systematic and has led to a rich collection of monotonicity formulas and quantitative consequences. In view of the intrinsic anisotropy of Finsler manifolds and the nonlinearity of the Finslerian Laplacian, extending the parabolic frequency method to the Finsler setting is significant not only for the study of nonlinear equations on metric measure spaces, but also for the development of the geometric analysis of Finsler manifolds.

From a technical point of view, frequency-function monotonicity and Carleman estimates represent two complementary approaches to unique continuation theory. The former is more geometric, whereas the latter is more analytic. Frequency monotonicity can be used to derive doubling inequalities, which in turn imply strong unique continuation properties \cite{GarofaloLin1986,GarofaloLin1987}. Carleman estimates, on the other hand, can be used to obtain three-spheres or three-circles inequalities, and hence unique continuation through propagation-of-smallness arguments \cite{MV2012}. However, the Carleman approach encounters substantial difficulties in the Finsler setting. First, the Finsler Laplacian is generally a nonlinear operator depending on a reference direction, which obstructs several standard linear arguments in Carleman estimates, including the expansion of conjugated operators, integration by parts, and the absorption of lower-order terms. Second, the fundamental tensor $g_{\nabla u}$ has no distinguished smooth value on the critical set $\left\{\nabla u=0\right\}$, so the global use of the linearized operator requires a weak reference-field convention. Moreover, the limited regularity of weak solutions and the singular behavior of distance functions or singular weight functions create further obstacles to a direct generalization. Therefore, studying nonlinear heat equations on Finsler manifolds through the frequency-function method is well adapted to the underlying geometric structure and may also provide useful insight for the development of Carleman estimates on Finsler metric measure spaces.

This generalization is also expected to contribute to the study of the geometry of Finsler metric measure spaces and of solutions to nonlinear equations on such spaces. Possible applications include generalizations of growth estimates for harmonic functions \cite{ColdingMinicozzi1997}, structural results for nodal sets of solutions to elliptic equations \cite{HardtSimon1989,Logunov2018}, blow-up analysis of singularities \cite{Huisken1990MCF,CM2012GenericMCF}, and quantitative stratification \cite{CheegerNaber2013, CheegerNaberValtorta2015}. Thus, in order to address these problems in Finsler geometric analysis and to better understand the global geometric properties of Finsler manifolds, a natural starting point is to develop an appropriate theory of parabolic frequency functions.

Historically, Almgren introduced the frequency function in his study of harmonic functions and regularity problems for minimal surfaces \cite{Almgren1977}. N. Garofalo and F.-H. Lin later developed this method for elliptic operators with variable coefficients and introduced a geometric variational approach to unique continuation \cite{GarofaloLin1986, GarofaloLin1987}. T. H. Colding and W. P. Minicozzi II applied frequency-function methods to harmonic functions on Riemannian manifolds with nonnegative Ricci curvature and Euclidean volume growth, obtaining control of polynomial growth, vanishing order, and Liouville-type properties \cite{ColdingMinicozzi1997}. C.-C. Poon was the first to extend the frequency function to the parabolic setting \cite{Poon1996Parabolic}. In Poon's parabolic frequency functional, the relevant $L^{2}$-norm is weighted by the backward heat kernel. Using the monotonicity of this parabolic frequency, Poon established unique continuation results and related quantitative estimates for the heat equation, and also discussed the relationship between the parabolic frequency method and Carleman methods in unique continuation theory. L. Ni further generalized parabolic frequency functions to heat equations on Riemannian and K\"ahler manifolds, proving monotonicity formulas by means of matrix Harnack estimates and Hamilton-type estimates \cite{Ni2015ParabolicFrequency}. More recently, N. Gigli and I. Y. Violo extended the Colding--Minicozzi monotonicity formula for positive harmonic functions to $\mathrm{RCD}(0,N)$ spaces and established rigidity and almost-rigidity results \cite{GigliViolo2023}. In this framework, rigidity means that equality in the monotonicity formula forces both the underlying space and the harmonic function to agree with the corresponding conical model, whereas almost rigidity asserts that an almost-equality case implies quantitative closeness to such a model.

Parabolic frequency functions also play an important role in the study of geometric flows. J. Baldauf, P. Ho, and T. Lee extended the notion of parabolic frequency to the setting of mean curvature flow. More precisely, they studied solutions to the heat equation along homothetically shrinking mean curvature flows induced by self-shrinkers and proved the monotonicity of the corresponding frequency function \cite{BalLee22}. X. Li and Q. S. Zhang established matrix Li--Yau--Hamilton-type estimates for the heat equation and its conjugate under Ricci flow, and applied these estimates to prove the monotonicity of parabolic frequency under a lower bound on sectional curvature and an upper bound on Ricci curvature \cite{LiZhang2024RicciFlow}. P. Xu and X.-A. Ren further introduced a modified frequency function in the setting of normalized Ricci flow on compact Riemannian manifolds with nonnegative sectional curvature and parallel Ricci curvature, as well as on closed surfaces with positive Euler characteristic, and proved analogous monotonicity and unique continuation results \cite{XuRen2025}.

In this paper, we establish monotonicity formulas for parabolic frequency functions and derive backward uniqueness on closed Finsler metric measure spaces. For a fixed reference measure, the frequency is nondecreasing along the heat equation when the metric is stationary, or when its evolution satisfies a suitable integral bound on the metric variation (see Theorem~\ref{thm1}). We also characterize equality and apply the result to Bao's normalized Finslerian Ricci flow (see Theorem~\ref{NRF1}). We then consider coupled metric measure flows, in which both the metric and the measure evolve, and refer to them as coupled flows below. Under an integral bound involving the metric variation and the change in weighted Ricci curvature, we prove monotonicity of a modified frequency for the heat equation
(see Theorem~\ref{thm:independent-heat}). For equations with a source controlled by the solution and its gradient, we obtain a differential inequality for the frequency and establish backward uniqueness under the corresponding assumptions. Thus a solution that vanishes at the final time must vanish throughout the preceding interval (see Section~\ref{sec:fixed-measure} and Theorem~\ref{thm:independent-backward}).


The situation becomes more delicate for coupled flows. In the Riemannian setting, weighted volume forms associated with the conjugate heat equation are used in parabolic frequency arguments \cite{JuliusKim2022}. Here we first construct local supersolutions and subsolutions of the conjugate heat equation from radial ordinary differential inequalities. We extend them to smooth positive global weights while preserving their respective differential inequalities. The time factors compensate only for the additional cutoff terms, rather than replacing the inequalities already established in the interior region. The resulting two densities need not be spatially proportional, and their comparison enters the frequency function. Under an upper bound on the integral of the metric-variation term and the difference of the weighted Ricci curvatures, we prove monotonicity for the heat equation and a differential inequality for equations with sources controlled by the solution and its gradient. To prove backward uniqueness, we first bound an energy quotient computed with the subsolution measure. This controls the term comparing the two weights up to the time at which the solution is assumed to vanish. The curvature assumption is imposed on the resulting global measures, and all global identities are used under the regularity assumptions stated below (see \eqref{eq:independent-regularity}).

The remainder of the paper is organized as follows. Section~\ref{sec:added-preliminaries} recalls the Finsler definitions and analytic identities used in the proofs. Section~\ref{sec:fixed-measure} establishes frequency monotonicity and backward uniqueness with a fixed reference measure, and applies these results to Bao's normalized Finslerian Ricci flow. Section~\ref{sec:added-coupled-examples} explains the role of an evolving measure and gives examples of coupled flows. Section~\ref{sec:local-conjugate-construction} constructs positive global supersolutions and subsolutions of the linearized conjugate heat equation. Finally, Section~\ref{sec:independent-frequency} uses these weights to prove frequency estimates and backward uniqueness under coupled flows.

\section{Preliminaries}\label{sec:added-preliminaries}

Throughout, $M$ is a connected closed smooth manifold of dimension $n$, and $d\mu$ is a smooth positive measure. For time-dependent metrics, gradients, norms, and curvature quantities are computed using $F(t)$ unless otherwise specified. We write $\tau=t_1-t>0$ and work on compact time intervals contained in $[t_0,t_1)$. We recall the definitions and identities needed below, following \cite{BaoChernShen2000,OhtaSturm2009HeatFlow,OhtaSturm2014Bochner}.

A Finsler metric $F:TM\to[0,\infty)$ is smooth outside the zero section, positively homogeneous of degree one in the tangent variable, and strongly convex. We write $v=v^i\partial_{x^i}\in T_xM$ for a tangent vector and use $(x^i,v^i)$ as local coordinates on $TM$. Its fundamental tensor is
\[
    g_{ij}(x,v)=\frac12\frac{\partial^2F^2}{\partial v^i\partial v^j}(x,v),\qquad v\ne0.
\]
Thus $g_v$ is an inner product on $T_xM$ and $g_v(v,v)=F^2(v)$. We write $F_v(w)=\sqrt{g_v(w,w)}$ for its induced norm. The dual norm is
\[
    F^*(x,\xi)=\sup_{F(x,v)=1}\xi(v).
\]
The gradient $\nabla u$ is determined by $du=g_{\nabla u}(\nabla u,\cdot)$ where $du\ne0$, and is defined to be zero where $du=0$. In particular,
\[
    du(\nabla u)=F^2(\nabla u)=F^{*2}(du).
\]
The map $du\mapsto\nabla u$ is generally nonlinear. Reversibility is not assumed, so $F(-v)$ need not equal $F(v)$.

In local coordinates, write $d\mu=\sigma(x)\,dx$, where $dx=dx^1\cdots dx^n$. The divergence and the Finsler Laplacian are defined by
\[
    \operatorname{div}_{\mu}X=\frac1\sigma\partial_{x^i}(\sigma X^i),\qquad
    \Delta_{\mu}u=\operatorname{div}_{\mu}(\nabla u).
\]
For a smooth test function $f$, integration by parts on the closed manifold gives
\[
    \int_M f\Delta_{\mu}u\,d\mu=-\int_Mdf(\nabla u)\,d\mu,
    \qquad
    \int_Mu\Delta_{\mu}u\,d\mu=-\int_MF^2(\nabla u)\,d\mu.
\]
The first identity also defines the weak Laplacian. We suppress the subscript $\mu$ when the reference measure is fixed and clear from the context. When the metric and measure depend on time, $\Delta_{\mu(t)}$ is formed using $F(t)$ and $\mu(t)$.

For a nonvanishing reference vector field $V$, the linearized gradient and Laplacian are
\[
    \nabla^V f=g^{ij}(x,V)\partial_{x^j}f\,\partial_{x^i},\qquad
    \Delta_{\mu}^V f=\operatorname{div}_{\mu}(\nabla^V f).
\]
Here $(g^{ij})$ is the inverse of $(g_{ij})$. These operators are linear in $f$ once $V$ is fixed. They satisfy
\[
    \int_M f\Delta_{\mu}^V k\,d\mu
    =-\int_M g_V(\nabla^V f,\nabla^V k)\,d\mu
    =\int_M k\Delta_{\mu}^V f\,d\mu.
\]
On $M_u=\{x\in M:du(x)\ne0\}$, we have $\nabla^{\nabla u}u=\nabla u$ and $\Delta_{\mu}^{\nabla u}u=\Delta_{\mu}u$. To define $\Delta_{\mu}^{\nabla u}$ on all of $M$, choose a measurable nonzero reference field that agrees with $\nabla u$ on $M_u$, as in \cite{OhtaSturm2014Bochner}, and use the associated divergence-form operator. Strong convexity and compactness give uniform ellipticity on compact smooth time intervals, although the coefficients need not be smooth across the critical set. All global identities involving this operator are understood weakly.

The measure is an additional part of the Finsler structure. Two standard choices are the Busemann--Hausdorff and Holmes--Thompson measures, given in local coordinates by
\[
\begin{aligned}
    d\mu_{\mathrm{BH},F}
    &=\frac{\omega_n}{\operatorname{vol}\{v\in T_xM:F(x,v)<1\}}\,dx,\\
    d\mu_{\mathrm{HT},F}
    &=\frac{\operatorname{vol}\{\xi\in T_x^*M:F^*(x,\xi)<1\}}{\omega_n}\,dx.
\end{aligned}
\]
Here the volumes in the numerators and denominators are Lebesgue volumes in the tangent coordinates and their dual coordinates, and $\omega_n$ is the Euclidean unit-ball volume. Both choices agree with the Riemannian volume in the Riemannian case. For a conformal change $F\mapsto e^{\phi(x)}F$, both measures are multiplied by $e^{n\phi(x)}$, as follows by rescaling the corresponding unit balls.

Changing the measure leaves $\nabla u$ unchanged, but changes its divergence. In particular, if $d\mathcal V=K\,d\mu$ for a smooth positive function $K$, then
\[
    \Delta_{\mathcal V}u=\Delta_{\mu}u+d\log K(\nabla u),\qquad
    \Delta_{\mathcal V}^V f=\Delta_{\mu}^V f+d\log K(\nabla^V f).
\]
We use these identities when the measure defining the frequency differs from the reference measure in the heat equation.

We next recall the weighted Ricci curvature. Fix $v\in T_xM\setminus\{0\}$ and let $\gamma$ be the $F$-geodesic with $\gamma(0)=x$ and $\dot\gamma(0)=v$. Along this geodesic, write
\[
    d\mu_{\gamma(s)}=e^{-\Psi(s)}d\operatorname{vol}_{g_{\dot\gamma(s)}}.
\]
The S-curvature and the weighted Ricci curvature used below are
\[
    \mathbf S_{\mu}(v)=\Psi'(0),\qquad
    \operatorname{Ric}_{\mu}^{\infty}(v)=\operatorname{Ric}(v)+\Psi''(0),
\]
where $\operatorname{Ric}(v)$ is the degree-two Finsler Ricci curvature. These quantities are defined to be zero at $v=0$. At a fixed time, replacing $d\mu$ by $d\mathcal V=K\,d\mu$ changes $\Psi$ to $\Psi-\log K\circ\gamma$. Consequently,
\[
    \operatorname{Ric}_{\mathcal V}^{\infty}(v)-\operatorname{Ric}_{\mu}^{\infty}(v)
    =-\left.\frac{d^2}{ds^2}\log K(\gamma(s))\right|_{s=0}.
\]
This formula expresses the curvature term in the estimates for coupled flows in terms of the spatial variation of $K$.

On $M_u$, let $\nabla^2u$ denote the Hessian formed with the Chern connection and reference direction $\nabla u$, and let $\|\nabla^2u\|_{HS(\nabla u)}$ be its Hilbert--Schmidt norm with respect to $g_{\nabla u}$. In the integrated identities below, the Hessian is extended by zero on $M\setminus M_u$. The Bochner--Weitzenb\"ock formula of Ohta--Sturm \cite{OhtaSturm2014Bochner} is
\[
    \Delta_{\mu}^{\nabla u}\left(\frac12F^2(\nabla u)\right)
    -d(\Delta_{\mu}u)(\nabla u)
    =\operatorname{Ric}_{\mu}^{\infty}(\nabla u)+\|\nabla^2u\|_{HS(\nabla u)}^2
\]
on $M_u$. Its integrated version applies, in particular, when $u\in H^2(M)\cap C^1(M)$ and $\Delta_{\mu}u\in H^1(M)$. Integrating over $M$ and using integration by parts gives
\[
    \int_M\|\nabla^2u\|_{HS(\nabla u)}^2\,d\mu
    =\int_M\left[(\Delta_{\mu}u)^2-\operatorname{Ric}_{\mu}^{\infty}(\nabla u)\right]d\mu.
\]
The same formula holds for any smooth positive measure $d\mathcal V$, with $\Delta_{\mathcal V}$ and $\operatorname{Ric}_{\mathcal V}^{\infty}$ in place of $\Delta_\mu$ and $\operatorname{Ric}_\mu^{\infty}$.

For a metric evolving by $\partial_t\log F=-H$ (see \eqref{F1}), the function $H=H(x,v,t)$ is homogeneous of degree zero in $v$ and may depend on its direction. Differentiating $\partial_tF^2=-2HF^2$ twice in the tangent variables yields
\[
    \partial_tg_{ij}=-2h_{ij},\qquad
    h_{ij}=\frac12\frac{\partial^2(HF^2)}{\partial v^i\partial v^j}.
\]
We write $h(v)=h_{ij}(x,v,t)v^iv^j$, with $h(0)=0$. Euler's identity gives $h(v)=H(x,v,t)F^2(v)$. When $F$ is stationary, $h=0$.

The following identity includes the time variation of both the metric and the gradient.

\begin{lemma}\label{lem1}
	If $F$ evolves by $\partial_t\log F=-H$ (see \eqref{F1}), then
	\[
		\partial_{t} F^{2}(\nabla u) = 2(du_{t}(\nabla u) + h(\nabla u)).
	\]
\end{lemma}
\begin{proof}
Differentiating $du=g_{\nabla u}(\nabla u,\cdot)$ and using the homogeneity identity $(\partial_{v^k}g_{ij})v^j=0$ gives, on $M_u$,
\[
    g_{\nabla u}(\partial_t\nabla u,\nabla u)
    =du_t(\nabla u)+2h(\nabla u).
\]
It follows that
\[
\begin{aligned}
    \partial_tF^2(\nabla u)
    &=-2h(\nabla u)+2g_{\nabla u}(\partial_t\nabla u,\nabla u)\\
    &=2du_t(\nabla u)+2h(\nabla u).
\end{aligned}
\]
\end{proof}

We assume sufficient regularity to differentiate the energies and apply the weak integration-by-parts identities. Precise regularity assumptions for the calculations along coupled flows are given below (see \eqref{eq:independent-regularity}). We first consider a fixed reference measure, for which Lemma~\ref{lem1} supplies the metric-variation term in the frequency calculation.

\section{Frequency estimates with a fixed measure}\label{sec:fixed-measure}

Throughout this section, the reference measure $d\mu$ is fixed, while the Finsler metric may depend on time. We begin with the heat equation.
\begin{theorem}\label{thm1}
	Let $(M,F,\mu)$ be a closed Finsler metric measure space, where $F$ is either time-independent or evolves by
	\begin{equation} \label{F1}
		\partial_{t} \log F = - H(t).
	\end{equation}
	Suppose that $u$ is a solution of the heat equation
	\begin{equation} \label{H1}
		\partial_{t} u = \Delta u,
	\end{equation}
	and assume that
	\begin{equation} \label{F1C}
		\int_{M} h(\nabla u) d \mu \leq \frac1{2\tau}\int_{M} F^{2}(\nabla u) d\mu,
	\end{equation}
	when $F$ evolves according to \eqref{F1}. Then the parabolic frequency
	\begin{equation} \label{U1}
		U(t) = - \tau \frac{\int_{M}F^{2}(\nabla u) d\mu}{\int_{M} u^2 d\mu}
	\end{equation}
	is nondecreasing on each interval on which $I(t)=\int_Mu^2\,d\mu>0$. For an evolving metric, equality $U'(t)=0$ holds if and only if
	\[
		\Delta u=c(t)u
	\]
	and equality holds in \eqref{F1C}. In this case,
	\[
		c(t)=\frac{U(t)}{\tau}.
	\]
	When $F$ is time-independent, $U'(t)=0$ if and only if $u$ is spatially constant. In particular, if $U'\equiv0$ on a time interval in the evolving case, then $U(t)=U_{0}$ and $c(t)=U_{0}/\tau$ on that interval.
\end{theorem}

\begin{proof}[Proof of Theorem \ref{thm1}]
	Set $I = \int_{M} u^2 d\mu$ and $D = -\tau \int_{M} F^{2}(\nabla u) d\mu \leq 0$, where $\tau = t_{1} - t > 0$. Differentiating $I$ and integrating by parts gives
	\begin{equation}
		I' = 2 \int_{M} u u_{t} d\mu = -2 \int_{M} F^{2}(\nabla u) d\mu = \frac{2}{\tau} D.
	\end{equation}
	Lemma~\ref{lem1} gives
	\begin{equation}
		\begin{split}
			D' &= \int_{M} F^{2}(\nabla u) d\mu - \tau \int_{M} \partial_{t} F^{2}(\nabla u) d\mu \\
			&= -\frac{D}{\tau} + 2\tau \int_{M} [u_{t}^{2} - h(\nabla u)] d\mu.
		\end{split}
	\end{equation}
	Since
	\[
		D=\tau\int_{M}u u_{t}\,d\mu,
	\]
	we obtain
	\begin{equation}\label{dU1}
		\begin{split}
			U'
			={}&\frac1{I^{2}} (D'I - I'D) \\
			={}&\frac{2\tau}{I^{2}}
			\left[
				\int_{M} u^2 d\mu \int_{M}u_{t}^{2} d\mu
				-
				\left(\int_{M} u u_{t} d\mu\right)^{2}
			\right]\\
			&+\frac1I
			\left[
				\int_{M}F^{2}(\nabla u)d\mu
				-
				2\tau\int_{M}h(\nabla u)d\mu
			\right]
			\geq 0.
		\end{split}
	\end{equation}
	The first term is nonnegative by the Cauchy--Schwarz inequality. The second is nonnegative by \eqref{F1C} when the metric evolves, and by $h=0$ when the metric is fixed. For an evolving metric, equality $U'(t)=0$ holds precisely when both inequalities are equalities, which means that
	\[
		u_t=\Delta u=c(t)u
	\]
	and
	\[
		\int_Mh(\nabla u)d\mu
		=
		\frac1{2\tau}\int_MF^2(\nabla u)d\mu.
	\]
	Substituting the eigenfunction relation into the expression for $U(t)$ gives
	\[
		U(t)=\tau c(t),
	\]
	and hence
	\begin{equation}
		c(t) = \frac{U(t)}{\tau}.
	\end{equation}
	If $F$ is time-independent, then $h=0$, and the second term in \eqref{dU1} is
	\[
		\frac1I\int_MF^2(\nabla u)d\mu.
	\]
	Thus $U'(t)=0$ if and only if $\nabla u\equiv0$, that is, $u$ is spatially constant.
\end{proof}

Monotonicity yields the following backward-uniqueness result.
\begin{corollary}
	Let the Finsler metric measure space $(M, F, \mu)$ and the function $u$ be as in Theorem \ref{thm1}. If the solution $u$ satisfies $u(\cdot, b) \equiv 0$ for some $b<t_{1}$, then $u \equiv 0$ for all $t \in [t_{0},b]$.
\end{corollary}
\begin{proof}
	Suppose, to the contrary, that $I(a)>0$ for some $a\in[t_{0},b)$. Since $I(b)=0$ and $I$ is continuous, let
	\[
		b_{0}:=\inf\{t\in(a,b]:I(t)=0\}.
	\]
	On $[a,b_{0})$, we have $I(t)>0$ and
	\begin{equation}
		(\log I)' = \frac{I'}{I} = \frac{2}{\tau}\frac{D}{I} = \frac{2}{\tau} U.
	\end{equation}
	For $s\in(a,b_{0})$, monotonicity gives
	\begin{equation}
		\begin{split}
			\log I(s) - \log I(a)
			&= 2 \int_{a}^{s} \frac{U(t)}{t_{1} - t} dt \\
			&\geq 2U(a)\int_{a}^{s} \frac{dt}{t_{1} - t} \\
			&= 2U(a) \log \frac{t_{1} - a}{t_{1} - s}.
		\end{split}
	\end{equation}
	Equivalently,
	\begin{equation}\label{ineq-IbIa-1}
		I(s) \geq I(a)\left(\frac{t_{1} - a}{t_{1} - s}\right)^{2U(a)} >0.
	\end{equation}
	Letting $s\rightarrow b_{0}$, the right-hand side of \eqref{ineq-IbIa-1} converges to a strictly positive number because $b_{0}\leq b<t_{1}$. This contradicts $I(b_{0})=0$. Hence $I(a)=0$ for every $a\in[t_{0},b]$, and therefore $u\equiv0$ on $M\times[t_{0},b]$.
\end{proof}

For a heat equation with a source $Q(u,\nabla u)$ (see \eqref{H2}), the frequency need not be monotone. Nevertheless, a differential inequality for the frequency still gives backward uniqueness when the source is bounded by the gradient and the solution.
\begin{theorem}
	Let $(M,F,\mu)$ be a closed Finsler metric measure space, where $F$ is either time-independent or satisfies \eqref{F1}. Suppose that $u$ is a sufficiently regular weak solution of
	\begin{equation} \label{H2}
		(\partial_{t} - \Delta) u = Q(u, \nabla u),
	\end{equation}
	on $M\times[t_{0},b]$, where $b<t_1$, so that the differentiations and integrations by parts below are justified, and assume that
	\begin{equation} \label{H2C}
		|Q(u,\nabla u)| \leq C(t)(F(\nabla u) + |u|),
	\end{equation}
	where $C(t)$ is bounded and nonnegative on $[t_{0},b]$.
	Assume also that \eqref{F1C} holds when the metric evolves. Then $u \equiv 0$ for all $t \in [t_{0}, b]$, provided $u(\cdot, b) \equiv 0$.
\end{theorem}

\begin{proof}
	Fix $a\in[t_{0},b)$. Suppose, to the contrary, that $I(a)>0$, and let
	\[
		b_{0}:=\inf\{t\in(a,b]:I(t)=0\}.
	\]
	Then $I(t)>0$ on $[a,b_{0})$, and the following computations hold on this interval.

	Differentiating $I$ and integrating by parts gives
	\begin{equation}
		\begin{split}
			I' = & 2\int_{M} u u_{t} d\mu \\
			= & -2\int_{M} F^{2}(\nabla u) d\mu + 2 \int_{M} u Q(u, \nabla u) d\mu \\
			\geq & \frac{2D}{\tau} - 2C(t) \int_{M} |u| (F(\nabla u) + |u|) d\mu \\
			\geq & \frac{2 + C(t)}{\tau} D - 3C(t)I,
		\end{split}
	\end{equation}
	where the last inequality follows from $2|u|F(\nabla u)\leq u^2+F^2(\nabla u)$. Therefore,
	\begin{equation}\label{ineq-logI'geq}
		(\log I)' \geq \frac{2 + C(t)}{\tau} U - 3C(t).
	\end{equation}
	We rewrite $D$ and $I'$ as
	\begin{equation}
		\begin{split}
			D = & \tau \int_{M} u \Delta u d\mu \\
			= & \tau \int_{M} u (u_{t}- Q) d\mu \\
			= & \tau \left(\int_{M} u(u_{t} - \frac{1}{2}Q) d\mu - \frac{1}{2}\int_{M} u Q d\mu\right),
		\end{split}
	\end{equation}
	and
	\begin{equation}
		\begin{split}
			I' = & -2\int_{M} F^{2}(\nabla u) d\mu + 2\int_{M} u Q d\mu \\
			= & 2\left(\int_{M} u(u_{t} - \frac{1}{2} Q) d\mu + \frac{1}{2}\int_{M} u Q d\mu\right),
		\end{split}
	\end{equation}
	so that
	\begin{equation}\label{DI'}
		DI' = 2\tau \left[\left(\int_{M} u(u_{t} - \frac{1}{2}Q)d\mu\right)^{2} - \frac{1}{4}\left(\int_{M} u Q d\mu\right)^{2}\right].
	\end{equation}
	For a time-independent metric, differentiation gives
	\begin{equation}
		\begin{split}
			D' = -\frac{D}{\tau} - 2\tau\int_{M}du_{t}(\nabla u) d\mu 
			= -\frac{D}{\tau} + 2\tau \int_{M}(u_{t} - Q)u_{t} d\mu.
		\end{split}
	\end{equation}
	For an evolving metric, Lemma~\ref{lem1} gives
	\begin{equation}
		\begin{split}
			D' = & -\frac{D}{\tau} - 2\tau \int_{M}[du_{t}(\nabla u) + h(\nabla u)] d\mu \\
			= & -\frac{D}{\tau} + 2\tau \int_{M}[(u_{t} - Q)u_{t} -h(\nabla u)] d\mu.
		\end{split}
	\end{equation}
	Using \eqref{F1C} in the evolving case, we obtain in both cases
	\begin{equation}\label{D'}
			D' \geq 2\tau \int_{M} (u_{t} - Q) u_{t} d\mu 
			=  2\tau \int_{M} \left[\left(u_{t} - \frac{1}{2}Q\right)^{2} - \frac{1}{4}Q^{2}\right] d\mu.
	\end{equation}
	It follows from \eqref{DI'} and \eqref{D'} that
	\begin{equation}\label{ineq-I2U'}
		\begin{split}
			I^{2}U' ={}& D'I - I'D \\
			\geq{}& 2\tau I\int_{M} (u_{t} - \frac{1}{2}Q)^{2} d\mu - \frac{\tau}{2} I \int_{M} Q^{2} d\mu\\
			& - 2\tau \left(\int_{M} u(u_{t} - \frac{1}{2}Q) d\mu\right)^{2} + \frac{\tau}{2}\left(\int_{M} uQ d\mu\right)^{2}.
		\end{split}
	\end{equation}
	The sum of the first and third terms on the right-hand side of \eqref{ineq-I2U'} is nonnegative by the Cauchy--Schwarz inequality. Therefore,
	\begin{equation}
		\begin{split}
			I^{2}U' \geq & -\frac{\tau}{2} I \int_{M} Q^{2} d\mu \\
			\geq & -\frac{C^{2}(t)\tau}{2} I \int_{M} (F(\nabla u) + |u|)^{2} d\mu \\
			\geq & C^{2}(t) I (D - \tau I).
		\end{split}
	\end{equation}
	Thus,
	\begin{equation}
		U' \geq C^{2}(t)(U - \tau) \geq C^{2}(t)(U - \tau(a)).
	\end{equation}
	Since $U\leq0$, we have $\tau(a)-U>0$, and hence
	\begin{equation}
		\frac{d}{dt}\log\bigl(\tau(a)-U(t)\bigr)
		=
		\frac{-U'(t)}{\tau(a)-U(t)}
		\leq C^{2}(t).
	\end{equation}
	Consequently,
	\begin{equation}
		\int_{a}^{t} C^{2}(s) ds \geq \log \frac{\tau(a) - U(t)}{\tau(a) - U(a)},
	\end{equation}
	which is equivalent to
	\begin{equation}\label{eq-U-lower-general}
		U(t) \geq \tau(a) - e^{\int_{a}^{t} C^{2}(s) ds} (\tau(a) - U(a)).
	\end{equation}

	Since $U(t)\leq0$, integrating \eqref{ineq-logI'geq} and using \eqref{eq-U-lower-general} gives, for every $s\in(a,b_0)$,
	\begin{equation}\label{ineq-I-general-final}
	\begin{aligned}
		\log\frac{I(s)}{I(a)}
		\geq{}&-3(s-a)\sup_{[a,b_0]}C+\left(2+\sup_{[a,b_0]}C\right)
		  \log\frac{t_1-a}{t_1-s}\times
		  \left[\tau(a)-e^{\int_a^{b_0}C^2(t)\,dt}
		  \bigl(\tau(a)-U(a)\bigr)\right].
	\end{aligned}
	\end{equation}
	The right-hand side has a finite limit as $s\rightarrow b_0$ because $C$ is bounded and $b_0\leq b<t_1$. Since $I(a)>0$, exponentiating \eqref{ineq-I-general-final} gives a strictly positive lower bound for $I(s)$ near $b_0$. This contradicts $I(b_0)=0$. Hence $I(a)=0$. Since $a\in[t_0,b)$ was arbitrary, $u\equiv0$ on $M\times[t_0,b]$.
\end{proof}

We next illustrate the fixed-measure results with Bao's normalized Finslerian Ricci flow. The hypothesis on the metric variation then becomes an integral condition on the Ricci curvature.

\begin{example}[Normalized Finslerian Ricci flow]
Bao introduced the normalized Finslerian Ricci flow in connection with the search for Einstein--Finsler metrics with constant Einstein factor \cite{Bao2005}. Write $F_t=F(t)$ and let $\mathcal R\mathrm{ic}_t(x,v)$ be its degree-zero Ricci scalar, so that
\[
    \mathcal R\mathrm{ic}_t(x,v)=\frac{\operatorname{Ric}_{ij}(x,v,t)v^iv^j}{F_t^2(x,v)}.
\]
The normalized flow is
\be\label{eq-NRF}
    \partial_t\log F_t=-\mathcal R\mathrm{ic}_t+\mathcal A(t),
\ee
where
\be\label{eq-NRF-C}
    \mathcal A(t)=\frac{1}{\operatorname{Vol}_{\sigma_t}(SM)}\int_{SM}\mathcal R\mathrm{ic}_t\,d\sigma_t,
    \qquad \operatorname{Vol}_{\sigma_t}(SM)=\int_{SM}d\sigma_t.
\ee
Here $SM=(TM\setminus\{0\})/\mathbb R_+$ is the sphere bundle, and $d\sigma_t$ is the volume form induced by the Sasaki-type metric
\[
    g_{ij}(x,v,t)\,dx^i\otimes dx^j
    +g_{ij}(x,v,t)\frac{\delta v^i}{F_t}\otimes\frac{\delta v^j}{F_t}.
\]
Thus $\mathcal A(t)$ is the average Ricci scalar over $SM$.

We use the Ricci tensor
\[
    \operatorname{Ric}_{ij}(x,v,t):=\frac12\frac{\partial^2}{\partial v^i\partial v^j}
       \left(F_t^2(x,v)\mathcal R\mathrm{ic}_t(x,v)\right).
\]
Since $F_t^2\mathcal R\mathrm{ic}_t$ is homogeneous of degree two,
\[
    \operatorname{Ric}_{ij}(x,v,t)v^iv^j=F_t^2(v)\mathcal R\mathrm{ic}_t(x,v)
\]
for every $v\in T_xM\setminus\{0\}$. Multiplying \eqref{eq-NRF} by $2F_t^2$ and differentiating twice in the fiber variables gives
\[
    \partial_tg_{ij}=-2\left(\operatorname{Ric}_{ij}-\mathcal A(t)g_{ij}\right).
\]
Hence the tensor $h_{ij}$ in $\partial_tg_{ij}=-2h_{ij}$ is
\[
    h_{ij}=\operatorname{Ric}_{ij}-\mathcal A(t)g_{ij}.
\]

To apply the fixed-measure results, we take the Holmes--Thompson measure of the initial metric,
\[
    d\mu:=d\mu_{\mathrm{HT},F_{t_0}},\qquad \partial_td\mu=0.
\]
All gradients and norms are computed using $F_t$. In particular,
\[
    \nabla u=\nabla^{F_t}u,\qquad
    \Delta_{\mu,t}u:=\operatorname{div}_{\mu}\bigl(\nabla^{F_t}u\bigr).
\]
Writing
\[
    \operatorname{Ric}_t(\nabla u):=\operatorname{Ric}_{ij}(x,\nabla u,t)(\nabla u)^i(\nabla u)^j,
\]
with value zero where $\nabla u=0$, we have
\[
    h_t(\nabla u)=\operatorname{Ric}_t(\nabla u)-\mathcal A(t)F_t^2(\nabla u).
\]
Condition \eqref{F1C} therefore takes the form
\be\label{eq-NRF-intRic}
    \int_M\left[\operatorname{Ric}_t(\nabla u)-\mathcal A(t)F_t^2(\nabla u)\right]d\mu
    \leq\frac{1}{2\tau(t)}\int_MF_t^2(\nabla u)\,d\mu,
\ee
where $\tau(t)=t_1-t$. Whenever $\int_MF_t^2(\nabla u)\,d\mu>0$, this is equivalent to
\[
    \frac{\displaystyle\int_M\operatorname{Ric}_t(\nabla u)\,d\mu}
         {\displaystyle\int_MF_t^2(\nabla u)\,d\mu}
    -\frac{\displaystyle\int_{SM}\mathcal R\mathrm{ic}_t\,d\sigma_t}
          {\operatorname{Vol}_{\sigma_t}(SM)}
    \leq\frac{1}{2\tau(t)}.
\]
The first term is the spatial average of the directional Ricci scalar, weighted by $F_t^2(\nabla u)$. The second is its average over $SM$. Thus \eqref{eq-NRF-intRic} bounds the difference between these two averages.

For a solution of $u_t=\Delta_{\mu,t}u$, the fixed-measure identities give
\[
    I'(t)=\frac{2}{\tau(t)}D(t)
\]
and
\[
    D'(t)=-\frac{D(t)}{\tau(t)}
    +2\tau(t)\int_M|\Delta_{\mu,t}u|^2\,d\mu
    -2\tau(t)\int_Mh_t(\nabla u)\,d\mu.
\]
Consequently,
\begin{align*}
    U'(t)={}&\frac{2\tau(t)}{I^2(t)}
    \left[I(t)\int_M|\Delta_{\mu,t}u|^2\,d\mu
    -\left(\int_Mu\Delta_{\mu,t}u\,d\mu\right)^2\right]\\
    &+\frac{1}{I(t)}\left[\int_MF_t^2(\nabla u)\,d\mu
    -2\tau(t)\int_Mh_t(\nabla u)\,d\mu\right].
\end{align*}
Under \eqref{eq-NRF-intRic}, the second term is nonnegative, and the Cauchy--Schwarz inequality gives
\[
    U'(t)\geq\frac{2\tau(t)}{I^2(t)}
    \left[I(t)\int_M|\Delta_{\mu,t}u|^2\,d\mu
    -\left(\int_Mu\Delta_{\mu,t}u\,d\mu\right)^2\right]\geq0.
\]
This proves the following specialization of Theorem~\ref{thm1}.

\begin{theorem}\label{NRF1}
Let $(M,F(t),\mu)$, $t\in[t_0,t_1)$, be a closed Finsler metric measure family satisfying the normalized Finslerian Ricci flow \eqref{eq-NRF}, with the fixed measure $d\mu=d\mu_{\mathrm{HT},F_{t_0}}$. Let $u\not\equiv0$ solve the heat equation \eqref{H1}. If \eqref{eq-NRF-intRic} holds along the flow, then the frequency $U(t)$ in \eqref{U1} is nondecreasing wherever $I(t)>0$.

At such a time, equality $U'(t)=0$ holds if and only if
\[
    \Delta_{\mu,t}u=c(t)u
\]
and
\[
    \int_Mh_t(\nabla u)\,d\mu=\frac{1}{2\tau(t)}\int_MF_t^2(\nabla u)\,d\mu.
\]
The second condition is equality in \eqref{eq-NRF-intRic}, and in this case $c(t)=U(t)/\tau(t)$.
\end{theorem}

The preceding argument also gives backward uniqueness.
\begin{corollary}\label{NRF2}
Let $(M,F(t),\mu)$, $t\in[t_0,t_1)$, be a closed Finsler metric measure family satisfying \eqref{eq-NRF}, with the fixed measure $d\mu=d\mu_{\mathrm{HT},F_{t_0}}$. Let
\[
    u:M\times[a,b]\longrightarrow\mathbb R,\qquad [a,b]\subset[t_0,t_1),
\]
solve the heat equation \eqref{H1}. If \eqref{eq-NRF-intRic} holds for every $t\in[a,b]$ and $u(\cdot,b)=0$, then $u\equiv0$ on $M\times[a,b]$.
\end{corollary}

The same fixed-measure estimates also give backward uniqueness in the presence of a source controlled by the solution and its gradient.

\begin{theorem}\label{NRF3}
Let $(M,F(t),\mu)$, $t\in[t_0,t_1)$, be a closed Finsler metric measure family satisfying \eqref{eq-NRF}, with the fixed measure $d\mu=d\mu_{\mathrm{HT},F_{t_0}}$. Let $u:M\times[a,b]\longrightarrow\mathbb R$, where $[a,b]\subset[t_0,t_1)$, solve \eqref{H2}. Assume that
\[
    u\in H^1\bigl([a,b];L_\mu^2(M)\bigr)\cap L^2\bigl([a,b];W_\mu^{2,2}(M)\bigr)
\]
and that a bounded nonnegative function $C:[a,b]\longrightarrow[0,\infty)$ satisfies
\[
    |Q(u,\nabla u)|\leq C(t)\bigl(F_t(\nabla u)+|u|\bigr)
\]
almost everywhere on $M\times[a,b]$. If \eqref{eq-NRF-intRic} holds for every $t\in[a,b]$ and $u(\cdot,b)=0$, then $u\equiv0$ on $M\times[a,b]$.
\end{theorem}

\end{example}

The conclusions above use the fixed measure $d\mu_{\mathrm{HT},F_{t_0}}$. Choosing the Holmes--Thompson or Busemann--Hausdorff measure of $F(t)$ instead allows the measure to vary with the metric, and its variation must then be included in the derivatives of $I$ and $D$. We therefore turn to evolving volume forms.

\section{Coupled flows}
\label{sec:added-coupled-examples}

\subsection{Evolving volume forms and examples of coupled flows}
In the preceding section, the Finsler metric evolves while the reference
measure remains fixed. This setting allows us to study the effect of the
metric variation, but it does not in general describe the evolution of the
metric together with its natural volume form. This distinction already
arises in Riemannian geometry. Under Ricci flow, the Riemannian volume form
evolves with the scalar curvature, while under a smooth mean curvature flow,
the induced volume form evolves with the squared norm of the mean curvature
vector \cite{Hamilton1982,Huisken1984}. Thus these flows, equipped with
their induced volume forms, cannot generally be treated within the
fixed-measure setting.

The same issue arises for a Finsler metric equipped with its
Busemann--Hausdorff or Holmes--Thompson measure. Both measures are determined
by the metric, so a change in the metric generally changes their densities
\cite[Sections~3.1 and~3.4]{AlvarezPaivaThompson2004}. Keeping the initial
measure fixed is a legitimate choice, but it is different from using the
corresponding measure of the evolving metric. To include both the classical
Riemannian flows and Finsler metric evolutions with their associated measures,
we therefore consider a coupled metric measure flow given by
\[
    \begin{cases}
        \partial_t\log F(t)=-H(t),\\
        \partial_t d\mu(t)=-\mathcal R(x,t)\,d\mu(t).
    \end{cases}
\]
We use this system in the construction of conjugate heat weights below
(see \eqref{F2}). Here \(H\) is homogeneous of degree zero in the tangent variable and may
depend on its direction, whereas \(\mathcal R\) is a function on
\(M\times[t_0,t_1)\). We do not impose a universal relation between these
two rates. The measure may be determined by \(F_t\), or it may include an
additional positive density depending on space and time. The fixed-measure
case is recovered when \(\mathcal R=0\).

We first describe the standard Riemannian and Finsler choices of measure,
and then give more explicit examples. Throughout this section, \(M\) is
closed and the evolutions are considered on an interval on which they
remain smooth.

\begin{example}[Riemannian metric and volume evolution]
\label{ex:added-riemannian}
Let \(g(t)\) satisfy \(\partial_tg=-2h\), and set
\[
    F_t(x,v)=\sqrt{g(t)_x(v,v)},
    \qquad
    d\mu_t=e^{-\psi(x,t)}\,d\operatorname{vol}_{g(t)}.
\]
Then
\[
    H(x,v,t)=\frac{h_x(v,v)}{g(t)_x(v,v)}.
\]
Since
\[
    \partial_t\sqrt{\det g(t)}
    =\frac12g^{ij}\partial_tg_{ij}\sqrt{\det g(t)}
    =-\operatorname{tr}_{g(t)}h\sqrt{\det g(t)},
\]
the measure evolves with rate
\[
    \mathcal R=\partial_t\psi+\operatorname{tr}_{g(t)}h.
\]
In particular, for Ricci flow with its Riemannian volume measure,
\[
    h=\operatorname{Ric}_{g(t)},
    \qquad
    \psi=0,
    \qquad
    \mathcal R=\operatorname{Scal}_{g(t)}.
\]
This is the usual volume evolution under Ricci flow
\cite{Hamilton1982}.

A smooth mean curvature flow of closed hypersurfaces gives another
example after the induced metric and volume form are pulled back to a
fixed manifold \(M\). Let
\(\iota_t:M\to\mathbb R^{n+1}\) satisfy
\(\partial_t\iota_t=\mathbf H\), where \(\mathbf H\) is the mean curvature
vector. Writing \(\mathrm{II}\) for the vector-valued second fundamental
form, with \(\mathbf H=\operatorname{tr}_{g(t)}\mathrm{II}\), we have
\[
    \partial_tg_{ij}
    =-2\langle\mathbf H,\mathrm{II}_{ij}\rangle,
    \qquad
    \partial_t d\operatorname{vol}_{g(t)}
    =-|\mathbf H|^2\,d\operatorname{vol}_{g(t)},
\]
as in \cite[Section~3]{Huisken1984}. Thus the induced metric and volume form
satisfy the coupled flow equations (see \eqref{F2}) with
\[
    h_{ij}=\langle\mathbf H,\mathrm{II}_{ij}\rangle,
    \qquad
    \mathcal R=|\mathbf H|^2.
\]
Here the mean curvature flow is an evolution of immersions, and it is
the induced metric--measure pair that enters our framework.
\end{example}

We next consider the measures determined by a Finsler metric. Unlike the
Riemannian case, the Busemann--Hausdorff and Holmes--Thompson choices need
not give the same measure or the same evolution rate.

\begin{example}[Finsler evolutions with their associated volume forms]
\label{ex:added-canonical-finsler}
Let \(F_t\) be a smooth Finsler metric evolution satisfying
\(\partial_t\log F_t=-H\). In fixed local coordinates, the standard volume
formulas are
\[
\begin{aligned}
    d\mu_{\mathrm{BH},t}
    &=\frac{\omega_n}
      {\operatorname{vol}\{v\in T_xM:F_t(x,v)<1\}}\,dx,\\
    d\mu_{\mathrm{HT},t}
    &=\frac{\operatorname{vol}\{\xi\in T_x^*M:F_t^*(x,\xi)<1\}}
      {\omega_n}\,dx,
\end{aligned}
\]
where \(\omega_n\) is the Euclidean unit ball volume and the fibre volumes
are computed in the corresponding tangent and cotangent coordinates
\cite[Sections~3.1 and~3.4]{AlvarezPaivaThompson2004}. Consequently,
\[
\begin{aligned}
    \mathcal R_{\mathrm{BH}}
    &=\partial_t\log\operatorname{vol}
      \{v\in T_xM:F_t(x,v)<1\},\\
    \mathcal R_{\mathrm{HT}}
    &=-\partial_t\log\operatorname{vol}
      \{\xi\in T_x^*M:F_t^*(x,\xi)<1\}.
\end{aligned}
\]
These are scalar functions on the base manifold. There is in general no
identity equating either rate with \(g_v^{ij}h_{ij}\) at an arbitrary
tangent direction, since this trace may depend on \(v\).

For the Busemann--Hausdorff measure, polar integration gives
\[
    \operatorname{vol}\{v:F_t(x,v)<1\}
    =\frac1n\int_{S^{n-1}}F_t(x,\theta)^{-n}\,d\theta,
\]
where \(S^{n-1}\) and \(d\theta\) are the Euclidean unit sphere and its
surface measure in the chosen tangent coordinates. Differentiating and
using \(\partial_t\log F_t=-H\), we obtain
\[
    \mathcal R_{\mathrm{BH}}(x,t)
    =n\frac{\displaystyle\int_{S^{n-1}}
       H(x,\theta,t)F_t(x,\theta)^{-n}\,d\theta}
      {\displaystyle\int_{S^{n-1}}F_t(x,\theta)^{-n}\,d\theta}.
\]
For Bao's normalized evolution \eqref{eq-NRF},
\(H=\mathcal R\mathrm{ic}_t-\mathcal A(t)\)
\cite{Bao2005}. The corresponding Busemann--Hausdorff rate is therefore
\[
    \mathcal R_{\mathrm{BH}}(x,t)
    =n\frac{\displaystyle\int_{S^{n-1}}
       \mathcal R\mathrm{ic}_t(x,\theta)F_t(x,\theta)^{-n}\,d\theta}
      {\displaystyle\int_{S^{n-1}}F_t(x,\theta)^{-n}\,d\theta}
      -n\mathcal A(t).
\]
The normalization term \(\mathcal A(t)\), being independent of \(x\),
does not force this rate to vanish at every point. More generally,
preservation of a total volume should not be confused with preservation
of the volume form pointwise.
\end{example}

The preceding examples use volume forms determined by the metric.
We now allow an additional positive density. A conformal evolution
provides a simple way to distinguish the change induced by the metric
from the change prescribed by that density.

\begin{example}[Conformal evolution with a time-dependent density]
\label{ex:added-conformal}
Let \((M,F_0,\mu_0)\) be a closed Finsler metric measure space, and let
\(\phi(x,t)\) and \(\psi(x,t)\) be smooth functions. Set
\[
    F_t(x,v)=e^{\phi(x,t)}F_0(x,v),
    \qquad
    d\mu_t=e^{n\phi(x,t)-\psi(x,t)}\,d\mu_0.
\]
The metric and measure equations hold with
\[
    H=-\partial_t\phi,
    \qquad
    h_{ij}=-(\partial_t\phi)g_{ij}(t),
    \qquad
    \mathcal R=\partial_t\psi-n\partial_t\phi.
\]
If \(\mu_0\) is the Busemann--Hausdorff or Holmes--Thompson measure of
\(F_0\), the corresponding measure of \(F_t\) is \(e^{n\phi}\mu_0\).
The factor \(e^{-\psi}\) therefore prescribes an additional density.
If \(F_0\) is non-Riemannian, then so is \(F_t\) at every time.

Since
\(\nabla^{F_t}u=e^{-2\phi}\nabla^{F_0}u\), the divergence formula gives
\[
    \Delta_{\mu_t,t}u
    =e^{-2\phi}\left[
      \Delta_{\mu_0,0}u+
      \bigl((n-2)d\phi-d\psi\bigr)(\nabla^{F_0}u)
      \right],
\]
where \(\Delta_{\mu_0,0}\) is the Laplacian associated with
\((F_0,\mu_0)\). In particular, taking \(\phi=0\) and
\(\psi(x,t)=(t-t_0)f(x)\), with \(f\) nonconstant, gives a stationary
metric with an evolving measure. In this case,
\[
    \mathcal R=f(x),
    \qquad
    \Delta_{\mu_t,t}u
    =\Delta_{\mu_0,0}u-(t-t_0)df(\nabla^{F_0}u).
\]
Thus the measure evolution alone introduces a time-dependent first-order
term into the heat equation.
\end{example}

For conformal changes, the Busemann--Hausdorff and Holmes--Thompson
measures acquire the same volume factor. The following Randers family
shows that their evolutions can differ under a nonconformal change of
the metric.

\begin{example}[A Randers evolution with two different volume forms]
\label{ex:added-randers}
Let \(M=\mathbb T^n=(\mathbb R/2\pi\mathbb Z)^n\), \(n\geq2\), and let
\(b(t)\) be smooth with \(|b(t)|<1\). Consider
\[
    F_t(x,v)=|v|+b(t)\cos x^1\,v^1.
\]
This defines a smooth strongly convex Randers metric. Wherever
\(b(t)\cos x^1\neq0\), it is nonreversible and hence non-Riemannian.
Its evolution is given by
\[
    H(x,v,t)=-\frac{b'(t)\cos x^1\,v^1}{F_t(x,v)},
    \qquad
    h(v)=-b'(t)\cos x^1\,v^1F_t(x,v).
\]
When \(b\) varies, the changing linear term gives a nonconformal
metric evolution.

For this family, the two standard measures are
\[
    d\mu_{\mathrm{BH},t}
    =\left(1-b(t)^2\cos^2x^1\right)^{(n+1)/2}\,dx,
    \qquad
    d\mu_{\mathrm{HT},t}=dx.
\]
To verify the first formula, fix \((x,t)\) and put
\(c=b(t)\cos x^1\). The tangent unit ball of \(|v|+cv^1\) is the
ellipsoid
\[
    (1-c^2)\left(v^1+\frac{c}{1-c^2}\right)^2
    +\sum_{j=2}^n(v^j)^2<\frac1{1-c^2},
\]
whose Euclidean volume is
\(\omega_n(1-c^2)^{-(n+1)/2}\). This gives the Busemann--Hausdorff
density. The dual unit ball is the Euclidean cotangent unit ball
translated by \(c\,dx^1\), because its support function is
\(|v|+cv^1\). Translation preserves its volume, which proves the
Holmes--Thompson formula.

Consequently,
\[
    \mathcal R_{\mathrm{BH}}(x,t)
    =\frac{(n+1)b(t)b'(t)\cos^2x^1}
      {1-b(t)^2\cos^2x^1},
    \qquad
    \mathcal R_{\mathrm{HT}}=0.
\]
The function \(\mathcal R_{\mathrm{BH}}\) is spatially nonconstant
whenever \(b(t)b'(t)\neq0\). Thus the same metric evolution belongs to
the fixed-measure setting when equipped with its Holmes--Thompson
measure, but gives an evolving measure when equipped with its
Busemann--Hausdorff measure. The associated heat operators and frequency
integrals therefore depend on the choice of measure as well as on the
metric evolution.
\end{example}

In the last example, the metric and the measure are determined by the
same positive function. We choose this function to solve the heat
equation on a fixed Finsler space, so that the construction gives an
explicit evolution of both parts of the metric--measure pair.

\begin{example}[Metric and volume evolution determined by a heat density]
\label{ex:added-heat-density}
On \(\mathbb T^n\), \(n\geq2\), fix \(\varepsilon>0\) and consider the
reversible Minkowski norm
\[
    F_0(v)=|v|+\varepsilon
      \sqrt{2(v^1)^2+\sum_{j=2}^n(v^j)^2}.
\]
This norm is smooth and strongly convex away from zero. Since the
cross term in \(F_0^2\) is not a quadratic form, \(F_0\) is
non-Riemannian. Its Busemann--Hausdorff measure \(d\mu_0\) is a constant
multiple of \(dx\).

For a smooth function \(f=f(x^1)\), symmetry and the Legendre relation
give
\[
    \nabla^{F_0}f
    =\frac{f_{x^1}}{(1+\sqrt2\,\varepsilon)^2}\,\partial_{x^1},
    \qquad
    \Delta_{\mu_0,0}f
    =\frac{f_{x^1x^1}}{(1+\sqrt2\,\varepsilon)^2}.
\]
It follows that
\[
    \rho(x,t)=1+\delta
    \exp\left(-\frac{t-t_0}{(1+\sqrt2\,\varepsilon)^2}\right)\cos x^1,
    \qquad 0<\delta<1,
\]
is a positive solution of
\(\partial_t\rho=\Delta_{\mu_0,0}\rho\) for \(t\geq t_0\).
Fix \(\alpha\neq0,1/n\) and define
\[
    F_t=\rho^\alpha F_0,
    \qquad
    d\mu_t=\rho\,d\mu_0.
\]
Then
\[
    H=-\alpha\frac{\partial_t\rho}{\rho},
    \qquad
    h_{ij}=-\alpha\frac{\partial_t\rho}{\rho}\,g_{ij}(t),
    \qquad
    \mathcal R=-\frac{\partial_t\rho}{\rho}.
\]
Both the metric and the volume form evolve, and \(\mathcal R\) is
spatially nonconstant. The total volume is nevertheless preserved,
since
\[
    \int_Md\mu_t
    =\int_M\rho\,d\mu_0
    =\int_Md\mu_0.
\]
The Busemann--Hausdorff measure of \(F_t\) is
\(\rho^{n\alpha}d\mu_0\). Thus the prescribed measure differs from it
by the spatially nonconstant factor \(\rho^{1-n\alpha}\).

The heat operator of the evolving pair is
\[
    \Delta_{\mu_t,t}u
    =\rho^{-2\alpha}\left[
      \Delta_{\mu_0,0}u+
      (1-2\alpha)d\log\rho(\nabla^{F_0}u)
      \right].
\]
The heat equation used to determine \(\rho\) is therefore the equation
on the fixed background \((F_0,\mu_0)\), rather than the heat equation
for the evolving pair \((F_t,\mu_t)\).
\end{example}

These examples show why the measure variation must be included in the
frequency calculation. All gradients and norms below are taken with
respect to \(F_t\). For a sufficiently regular solution of
\(u_t=\Delta_{\mu(t)}u\), differentiation and integration by parts give
\[
\begin{aligned}
    \frac{d}{dt}\int_Mu^2\,d\mu(t)
    ={}&-2\int_MF^2(\nabla u)\,d\mu(t)
       -\int_M\mathcal R u^2\,d\mu(t).
\end{aligned}
\]
Similarly, Lemma~\ref{lem1} gives
\[
\begin{aligned}
    \frac{d}{dt}\int_MF^2(\nabla u)\,d\mu(t)
    ={}&-2\int_M(\Delta_{\mu(t)}u)^2\,d\mu(t)
       +2\int_Mh(\nabla u)\,d\mu(t)\\
       &-\int_M\mathcal R F^2(\nabla u)\,d\mu(t).
\end{aligned}
\]
When \(\mathcal R\) depends only on time, the terms due to the measure
variation cancel in the derivative of the quotient
\[
    \frac{\displaystyle\int_MF^2(\nabla u)\,d\mu(t)}
      {\displaystyle\int_Mu^2\,d\mu(t)},
\]
whenever the denominator is nonzero. For spatially varying
\(\mathcal R\), the two measure terms have different weights, so this
cancellation generally fails. Preservation of total volume alone does
not remove them.

The change of measure also affects the heat operator. Writing
\(d\mu_t=e^{-\psi(x,t)}d\mu_0\) relative to a fixed smooth measure,
we have
\[
    \Delta_{\mu(t)}u
    =\Delta_{\mu_0,t}u-d\psi(\nabla u),
\]
where \(\Delta_{\mu_0,t}\) uses the evolving metric \(F_t\) and the fixed
measure \(\mu_0\). Hence
\[
    u_t-\Delta_{\mu(t)}u=0
    \quad\Longleftrightarrow\quad
    u_t-\Delta_{\mu_0,t}u=-d\psi(\nabla u).
\]
Thus replacing the evolving measure by a fixed reference measure
introduces a first-order term, rather than recovering the fixed-measure
heat equation without modification. The energy integrals and the
weighted Ricci curvature entering the frequency estimates also depend
on the chosen measure.

To handle the additional terms in the energy identities, we work with
a smooth positive weight \(K\). Integration by parts with the linearized
Laplacian gives
\[
\begin{aligned}
    \frac{d}{dt}\int_Mu^2K\,d\mu(t)
    ={}&-2\int_MF^2(\nabla u)K\,d\mu(t)\\
       &+\int_Mu^2\left(
         \partial_tK+\Delta_{\mu(t)}^{\nabla u}K-\mathcal R K
         \right)d\mu(t).
\end{aligned}
\]
The second integral vanishes when \(K\) satisfies the linearized conjugate
heat equation \(\partial_tK+\Delta_{\mu(t)}^{\nabla u}K-\mathcal R K=0\)
(see \eqref{CH}). If \(K\) is
a supersolution or a subsolution, its sign instead gives a one-sided
energy estimate. Before constructing these weights, we record existence for the associated heat equations.

\subsection{Existence and local smoothness of heat solutions under coupled flows}\label{subsec:coupled-heat-existence}

Along the prescribed smooth coupled flow, the nonlinear heat equation admits a variational solution. We use the time-dependent convex gradient-flow theorem of Kopfer \cite[Theorem~5.5]{Kopfer2018GradientFlow} in the scalar space $L^2(M,\mu(t_0))$, following the energy construction of Ohta--Sturm in the static setting \cite[Theorem~3.4]{OhtaSturm2009HeatFlow}. The full verification and the implicit time construction are given in Appendix~\ref{app:heat-existence}. We then establish local smooth existence near a point where the initial differential is nonzero.

\begin{prop}[Existence along a smooth coupled flow]
\label{prop:coupled-heat-existence}
Let $(M,F(t),\mu(t))$, $t\in[t_0,b]$ with $t_0<b<t_1$, be a prescribed smooth family of closed Finsler metric measure spaces satisfying \eqref{F2}. Assume that each $F(t)$ is strongly convex and each $\mu(t)$ is a smooth positive measure. For every $u_0\in H^1(M)$, the problem
\begin{equation}\label{eq:coupled-heat-cauchy}
    u_t=\Delta_{\mu(t)}u,\qquad u(\cdot,t_0)=u_0
\end{equation}
has a unique weak solution in the class
\begin{equation}\label{eq:coupled-heat-existence-regularity}
    u\in H^1((t_0,b);L^2(M,\mu(t_0)))
       \cap L^\infty((t_0,b);H^1(M)).
\end{equation}
In particular, $u\in C([t_0,b];L^2(M,\mu(t_0)))$, and the initial condition holds in this space. The equation is understood through
\[
    \int_M u_t\varphi\,d\mu(t)
    =-\int_M d\varphi(\nabla u)\,d\mu(t),
    \qquad \varphi\in H^1(M),
\]
for almost every time. Nonnegative initial data give a nonnegative solution.
\end{prop}

\begin{proof}
See Appendix~\ref{app:heat-existence}, where the hypotheses of Kopfer's theorem are checked for the Finsler energy and its subgradient is identified with $-\Delta_{\mu(t)}$. The implicit approximation also gives the stated regularity and preservation of nonnegativity.
\end{proof}

This variational construction does not assert smoothness at critical points. The corresponding distinction is already present for a fixed Finsler metric \cite[Example~4.3 and Remark~4.10]{OhtaSturm2009HeatFlow}. A local classical solution can be constructed when the initial differential is nonzero.

\begin{prop}[Local smooth solutions]\label{prop:local-smooth-heat}
Under the assumptions of Proposition~\ref{prop:coupled-heat-existence}, let $u_0$ be smooth near $x_0\in M$ and suppose that $du_0(x_0)\ne0$. There are a neighborhood $V$ of $x_0$ and a number $\delta>0$, with $t_0+\delta\leq b$, such that
\[
    u_t=\Delta_{\mu(t)}u,\qquad u(\cdot,t_0)=u_0
\]
has a smooth solution on $V\times[t_0,t_0+\delta]$ satisfying $du\ne0$. If $u_0(x_0)>0$, the neighborhood and time interval may be chosen so that $u>0$ there.
\end{prop}

\begin{proof}
In a coordinate neighborhood write $d\mu(t)=\sigma(x,t)\,dx$. For a nonzero covector $\xi$, set
\[
\begin{aligned}
    a^{ij}(x,t,\xi)
    &=\frac12\frac{\partial^2 F_t^{*2}}{\partial\xi_i\partial\xi_j}(x,\xi),\\
    B(x,t,\xi)
    &=\frac1{2\sigma(x,t)}
       \frac{\partial}{\partial x^i}
       \left(\sigma(x,t)
             \frac{\partial F_t^{*2}}{\partial\xi_i}(x,\xi)\right),
\end{aligned}
\]
where $\xi$ is held fixed in the $x$-derivative. Expanding the divergence gives
\[
    \Delta_{\mu(t)}u
    =a^{ij}(x,t,Du)\,\partial_{x^i x^j}u+B(x,t,Du)
\]
wherever $du\ne0$. The coefficients are smooth for $\xi\ne0$, and $a^{ij}$ is uniformly positive definite on compact base and time regions. This follows from strong convexity and the zero homogeneity of $a^{ij}$ in $\xi$.

Choose $V$ with compact closure in the coordinate neighborhood so that $du_0$ is nonzero near $\overline V$, and extend $u_0$ smoothly to a function $\widetilde u_0$ on $M$ that agrees with it there. Fix an auxiliary smooth Riemannian metric $\widehat g$. Choose a smooth function $0\leq\eta\leq1$ on $T^*M$, compactly supported over the coordinate neighborhood and away from the zero section, which equals one on a neighborhood of
\[
    \{(x,du_0(x))\mid x\in\overline V\}.
\]
Consider on the closed manifold the modified equation
\[
    \partial_t\widetilde u
    =\eta(x,d\widetilde u)\Delta_{\mu(t)}\widetilde u
     +(1-\eta(x,d\widetilde u))\Delta_{\widehat g}\widetilde u,
    \qquad \widetilde u(\cdot,t_0)=\widetilde u_0.
\]
The first term is defined to be zero near the zero section, where $\eta$ vanishes. The expanded equation has smooth coefficients everywhere, and its principal coefficient is
\[
    \eta a^{ij}+(1-\eta)\widehat g^{ij}.
\]
This is uniformly positive definite. The short-time existence theorem for smooth quasilinear parabolic equations on closed manifolds \cite[Theorem~1.1]{MantegazzaMartinazzi2012} therefore gives a smooth solution. Since $\widetilde u(\cdot,t)\to\widetilde u_0$ in $C^1(M)$ as $t\rightarrow t_0$, after shortening the time interval we have $\eta(x,d\widetilde u)=1$ and $d\widetilde u\ne0$ on $V$. The restriction of $\widetilde u$ is the required local solution of the original equation. Positivity follows by continuity after shrinking $V$ and the time interval.
\end{proof}

For example, in coordinates centered at $x_0$ and on a sufficiently small neighborhood, the initial function $u_0=1+x^1$ is positive and has nonzero differential. Proposition~\ref{prop:local-smooth-heat} then gives a positive nonconstant local smooth solution along the prescribed flow.

Proposition~\ref{prop:coupled-heat-existence} gives a global weak solution, whereas Proposition~\ref{prop:local-smooth-heat} concerns a local smooth solution. For the global frequency identities we retain the additional regularity in \eqref{eq:independent-regularity} and the smooth noncritical neighborhood of the cutoff support. The weights used below also have spatial derivatives supported in such a region and satisfy global differential inequalities. The radial barriers and cutoff construction in Section~\ref{sec:local-conjugate-construction} provide these properties.

\section{Conjugate heat equation and new volume forms}\label{sec:local-conjugate-construction}

We now study the conjugate heat equation associated with the coupled flow
\begin{equation}\label{F2}
    \begin{cases}
        \partial_t\log F(t)=-H(t),\\
        \partial_t d\mu(t)=-\mathcal R(x,t)\,d\mu(t).
    \end{cases}
\end{equation}
Throughout this section, $M$ is closed, and the metric and measure evolve smoothly on $[t_0,b]\subset[t_0,t_1)$. As above, $H$ may depend on the tangent direction, while $\mathcal R$ depends only on $(x,t)$. We use the tensor $h$ defined in the preliminaries. For a given solution $u$, the associated conjugate equation is
\begin{equation}\label{CH}
    \partial_t f=-\Delta_{\mu(t)}^{\nabla u}f+\mathcal R f.
\end{equation}
We first construct positive local super- and subsolutions $k_1,k_2$ from differential inequalities in the distance variable. A supersolution satisfies $\partial_t f+\Delta_{\mu(t)}^{\nabla u}f-\mathcal Rf\geq0$, and a subsolution satisfies the reverse inequality. We then use a cutoff to extend these functions to the constant density $1$ outside a regular region. Multiplication by suitable exponential factors in time controls the resulting error terms and gives the positive global weights $K_1,K_2$.

Fix $p\in M$ and put $r^+(x)=d^+_{F(t_0)}(p,x)$. Let
\[
    M_p^{\mathrm{reg}}
    =M\setminus\bigl(\{p\}\cup\operatorname{Cut}_{F(t_0)}^+(p)\bigr).
\]
Choose a relatively open space-time region
$\Omega_0\subset M\times[t_0,b]$ such that $u$ is smooth, $du\ne0$,
and $x\in M_p^{\mathrm{reg}}$ on a neighborhood of its closure in
$M\times[t_0,b]$. All estimates below are taken on
$\overline{\Omega_0}$, and the cutoff will be supported in a smaller
region. Define
\[
    \mathfrak{R}(x,t)=\int_{t_0}^t\mathcal R(x,s)\,ds.
\]
Strong convexity and compactness imply
\begin{equation}\label{positive-distance-coefficient}
    \min_{\overline{\Omega_0}}
    F_{\nabla u}^2(\nabla^{\nabla u}r^+)>0.
\end{equation}
Indeed, $F_{t_0}^*(dr^+)=1$ where the distance is smooth, and the
dual fundamental tensors are uniformly equivalent to a fixed positive
definite metric on the chosen compact region. Choose constants
\begin{align}
    L_1&\geq\sup_{\overline{\Omega_0}}
    \frac{e^{\mathfrak{R}}
          \bigl(|\Delta_{\mu(t)}^{\nabla u}\mathfrak{R}|
          +F_{\nabla u}^2(\nabla^{\nabla u}\mathfrak{R})\bigr)}
         {F_{\nabla u}^2(\nabla^{\nabla u}r^+)},
    \qquad L_1\geq0,\label{L1}\\
    L_2&\geq\sup_{\overline{\Omega_0}}
    \frac{|\mathcal R|}
         {F_{\nabla u}^2(\nabla^{\nabla u}r^+)},
    \qquad L_2>0,\label{L2}\\
    L_3&\geq\sup_{\overline{\Omega_0}}
    \frac{\max\{0,-\Delta_{\mu(t)}^{\nabla u}r^+\}}
         {F_{\nabla u}^2(\nabla^{\nabla u}r^+)},
    \qquad L_3>0.\label{L3}
\end{align}
All derivatives in \eqref{L1}--\eqref{L3} are computed with reference
field $\nabla u$. Smoothness and
\eqref{positive-distance-coefficient} ensure that the displayed
suprema are finite. After choosing $L_2$, increase $L_3$ if necessary
so that $L_2\leq L_3^2/4$.

With these constants fixed, restrict to a smaller region whose distance values lie in $[r_0,r_0+\ell]$, where $r_0>0$, and write $s=r^+-r_0\in[0,\ell]$. Since $ds=dr^+$ and $\Delta_{\mu(t)}^{\nabla u}s=\Delta_{\mu(t)}^{\nabla u}r^+$, the radial differential inequalities are unchanged. By shrinking the region, we may choose $\ell$ to satisfy the conditions below.

\subsection{A local supersolution}

We seek a supersolution of the form
\begin{equation}\label{K1}
    k_1(x,t)=e^{\mathfrak{R}(x,t)}+\mathcal C_1(s(x)).
\end{equation}
Since $\partial_t e^{\mathfrak{R}}=\mathcal R e^{\mathfrak{R}}$, the chain rule gives
\begin{equation}\label{orderequation}
\begin{aligned}
    &\partial_tk_1+\Delta_{\mu(t)}^{\nabla u}k_1-\mathcal Rk_1\\
    &\quad=\mathcal C_1^{\prime\prime}(s)F_{\nabla u}^2(\nabla^{\nabla u}r^+)
       +\mathcal C_1^{\prime}(s)\Delta_{\mu(t)}^{\nabla u}r^+
       -\mathcal R\mathcal C_1(s)\\
    &\qquad+e^{\mathfrak{R}}\left[
       \Delta_{\mu(t)}^{\nabla u}\mathfrak{R}
       +F_{\nabla u}^2(\nabla^{\nabla u}\mathfrak{R})\right].
\end{aligned}
\end{equation}
Primes denote derivatives with respect to $s$. If $\mathcal C_1>0$ and $\mathcal C_1^{\prime}\geq0$, the bounds \eqref{L1}--\eqref{L3} give
\[
\begin{aligned}
    \partial_tk_1+\Delta_{\mu(t)}^{\nabla u}k_1-\mathcal Rk_1
    \geq F_{\nabla u}^2(\nabla^{\nabla u}r^+)
    \left(\mathcal C_1^{\prime\prime}-L_3\mathcal C_1^{\prime}-L_2\mathcal C_1-L_1\right).
\end{aligned}
\]
It is therefore enough to solve
\begin{equation}\label{eq-super}
    \frac{d^2f}{ds^2}-L_3\frac{df}{ds}-L_2f-L_1=0.
\end{equation}
The general solution is
\begin{equation}\label{C1}
    \mathcal C_1(s)=\alpha_1e^{\lambda_1s}+\alpha_2e^{\lambda_2s}-\frac{L_1}{L_2},
\end{equation}
where
\[
    \lambda_1=\frac{L_3+\sqrt{L_3^2+4L_2}}2>0,
    \qquad
    \lambda_2=\frac{L_3-\sqrt{L_3^2+4L_2}}2<0.
\]
Choose
\[
    \alpha_2\geq0,\qquad
    \alpha_1>\max\left\{-\frac{\lambda_2}{\lambda_1}\alpha_2,\frac{L_1}{L_2}\right\}.
\]
Then $\mathcal C_1(0)>0$, and
\[
    \mathcal C_1^{\prime}(s)
    =\alpha_1\lambda_1e^{\lambda_1s}+\alpha_2\lambda_2e^{\lambda_2s}\geq0
    \qquad(s\geq0).
\]
The positive term is increasing, while the magnitude of the negative term is decreasing. Hence $\mathcal C_1>0$ on $[0,\ell]$ and
\begin{equation}\label{K1C}
    \partial_tk_1+\Delta_{\mu(t)}^{\nabla u}k_1-\mathcal Rk_1\geq0
    \qquad\text{on the chosen regular region}.
\end{equation}
The cutoff used below is supported in this regular region.

\subsection{A local subsolution}

Similarly, we seek a subsolution of the form
\begin{equation}\label{subsolution-ansatz}
    k_2(x,t)=e^{\mathfrak{R}(x,t)}+\mathcal C_2(s(x)).
\end{equation}
The required inequality is
\begin{equation}\label{subsolution-pointwise-condition}
\begin{aligned}
    &\partial_tk_2+\Delta_{\mu(t)}^{\nabla u}k_2-\mathcal Rk_2\\
    &\quad=\mathcal C_2^{\prime\prime}(s)F_{\nabla u}^2(\nabla^{\nabla u}r^+)
      +\mathcal C_2^{\prime}(s)\Delta_{\mu(t)}^{\nabla u}r^+
      -\mathcal R\mathcal C_2(s)\\
    &\qquad+e^{\mathfrak{R}}\left[\Delta_{\mu(t)}^{\nabla u}\mathfrak{R}
       +F_{\nabla u}^2(\nabla^{\nabla u}\mathfrak{R})\right]\leq0.
\end{aligned}
\end{equation}
For $\mathcal C_2>0$ and $\mathcal C_2^{\prime}\leq0$, multiplying the lower bound for $\Delta_{\mu(t)}^{\nabla u}r^+$ by $\mathcal C_2^{\prime}$ reverses the inequality. Hence
\[
\begin{aligned}
    \partial_tk_2+\Delta_{\mu(t)}^{\nabla u}k_2-\mathcal Rk_2
    \leq F_{\nabla u}^2(\nabla^{\nabla u}r^+)
       \left(\mathcal C_2^{\prime\prime}-L_3\mathcal C_2^{\prime}+L_2\mathcal C_2+L_1\right).
\end{aligned}
\]
It suffices to find a function satisfying
\begin{equation}\label{subsolution-control-ode}
    \frac{d^2f}{ds^2}-L_3\frac{df}{ds}+L_2f+L_1\leq0,
    \qquad f>0,\quad \frac{df}{ds}\leq0.
\end{equation}
Recall that
\begin{equation}\label{subsolution-discriminant}
    L_2\leq\frac{L_3^2}{4}.
\end{equation}
On an interval with $\ell<2/L_3$, choose
\[
    c>
    \frac{4L_1}{L_3^2}\frac{e^{-L_3\ell/2}}{2/L_3-\ell}
\]
and set
\begin{equation}\label{critical-subsolution}
    \mathcal C_2(s)=c\left(\frac2{L_3}-s\right)e^{L_3s/2}
       -\frac{4L_1}{L_3^2}.
\end{equation}
Direct differentiation gives
\[
    \mathcal C_2^{\prime}(s)=-\frac{cL_3}{2}s e^{L_3s/2}\leq0,
    \qquad
    \mathcal C_2^{\prime\prime}-L_3\mathcal C_2^{\prime}
         +\frac{L_3^2}{4}\mathcal C_2+L_1=0.
\]
The choice of $c$ gives $\mathcal C_2(\ell)>0$. Since $\mathcal C_2$ is nonincreasing, it is positive throughout $[0,\ell]$, and
\[
    \mathcal C_2^{\prime\prime}-L_3\mathcal C_2^{\prime}+L_2\mathcal C_2+L_1
      =-\left(\frac{L_3^2}{4}-L_2\right)\mathcal C_2\leq0.
\]
This proves \eqref{subsolution-pointwise-condition} whenever $\ell<2/L_3$.

A second family of local subsolutions is obtained by joining an exponential expression to a trigonometric one. Choose
\[
    0<\varepsilon<\frac{2}{L_{3}},
    \qquad
    s_{0}:=\frac{2}{L_{3}}-\varepsilon,
    \qquad
    s_{0}<\ell,
\]
and take $\varepsilon_{1},\varepsilon_{2}>0$ such that
\[
    \varepsilon_{1}+\varepsilon_{2}<\frac{\pi}{3}.
\]
Choose $\beta>0$ satisfying
\begin{equation}\label{restriction-for-beta}
    \frac{L_3}{2}\cot\left(\frac{\frac{\pi}{3}-\varepsilon_1-\varepsilon_2}{2}\right)
    \leq\beta\leq
    \frac{\frac{\pi}{3}-\varepsilon_1-\varepsilon_2}{2(\ell-s_0)}.
\end{equation}
The interval in \eqref{restriction-for-beta} is nonempty whenever
\begin{equation}\label{subsolution-diameter-condition}
    \ell-s_0\leq
    \frac{\frac{\pi}{3}-\varepsilon_1-\varepsilon_2}{L_3}
    \tan\left(\frac{\frac{\pi}{3}-\varepsilon_1-\varepsilon_2}{2}\right).
\end{equation}
For such a choice of $\beta$, \eqref{restriction-for-beta} gives
\[
    \beta(\ell-s_{0}) \leq \frac{
        \frac{\pi}{3}-\varepsilon_{1}-\varepsilon_{2}
    }{2}
\]
and
\[
    \arctan\frac{L_{3}}{2\beta} \leq \frac{
        \frac{\pi}{3}-\varepsilon_{1}-\varepsilon_{2}
    }{2}.
\]
Consequently,
\[
    \beta(\ell-s_{0}) + \arctan\frac{L_{3}}{2\beta} \leq \frac{\pi}{3}-\varepsilon_{1}-\varepsilon_{2}.
\]
These inequalities allow a phase for which the trigonometric expression is positive and decreasing on $[s_0,\ell]$.

To ensure that the function on the first interval is also decreasing, impose the additional condition
\begin{equation}\label{transmission-compatibility}
    \beta(\ell-s_0)+\max\left\{
        \arctan\frac{L_3}{2\beta},\arctan\frac{1}{\beta\varepsilon}
    \right\}<\frac{\pi}{3}-\varepsilon_1.
\end{equation}
The two conditions are compatible when $\ell-s_0$ is sufficiently small. First choose $\beta$ large enough to satisfy the lower bound in \eqref{restriction-for-beta} and to make the maximum in \eqref{transmission-compatibility} smaller than $\pi/3-\varepsilon_1$. Then choose $\ell-s_0$ small enough to satisfy both upper bounds. Since $\pi/3-\varepsilon_1<\pi/2$, condition \eqref{transmission-compatibility} also implies
\[
    \beta(\ell-s_0)<\arctan(\beta\varepsilon)<\beta\varepsilon,
    \qquad\text{and hence}\qquad \ell<\frac2{L_3}.
\]
Thus this construction also requires $\ell<2/L_3$. Under \eqref{transmission-compatibility}, choose $\phi$ such that
\begin{equation}\label{subsolution-phase-choice}
    \beta\ell-\left(\frac{\pi}{3}-\varepsilon_1\right)
    <\phi<\beta s_0-\max\left\{
        \arctan\frac{L_3}{2\beta},\arctan\frac{1}{\beta\varepsilon}
    \right\}.
\end{equation}
Consequently,
\[
    \arctan\frac{L_{3}}{2\beta} < \beta s-\phi < \frac{\pi}{3}-\varepsilon_{1}
    \qquad
    \text{for all }s\in[s_{0},\ell],
\]
and
\[
    \tan(\beta s_{0}-\phi) > \frac{1}{\beta\varepsilon}.
\]

Choose $\rho>0$ such that
\begin{equation}\label{subsolution-rho-choice}
    \rho>\frac{4L_1}{L_3^2}
    \frac{e^{-L_3\ell/2}}{\cos(\beta\ell-\phi)},
\end{equation}
and set
\begin{equation}\label{complete-subsolution-profile}
    \mathcal C_2(s)=
    \begin{cases}
        \begin{aligned}
        &\rho e^{L_3s/2}\bigl[\cos(\beta s_0-\phi)\\
        &\qquad+\beta(s_0-s)\sin(\beta s_0-\phi)\bigr]
        -\dfrac{4L_1}{L_3^2},
        \end{aligned}
        &0\leq s\leq s_0,\\[0.6em]
        \rho e^{L_3s/2}\cos(\beta s-\phi)-\dfrac{4L_1}{L_3^2},
        &s_0\leq s\leq\ell.
    \end{cases}
\end{equation}
Direct differentiation gives
\[
    \mathcal C_{2}(s_{0}^{-}) = \mathcal C_{2}(s_{0}^{+}),
    \qquad
    \mathcal C_{2}'(s_{0}^{-}) = \mathcal C_{2}'(s_{0}^{+}).
\]
Hence $\mathcal C_2$ is $C^1$ across $s=s_0$.

For $s\in[s_{0},\ell]$,
\[
\begin{aligned}
    \mathcal C_{2}'(s)
    ={}
    \rho e^{\frac{L_{3}}{2}s}
    \left[
        \frac{L_{3}}{2}\cos(\beta s-\phi) - \beta\sin(\beta s-\phi)
    \right]
    <0.
\end{aligned}
\]
The choice \eqref{subsolution-rho-choice} gives
$\mathcal C_{2}(\ell)>0$, and hence
\[
    \mathcal C_{2}(s)>0
    \qquad
    \text{for every }s\in[s_{0},\ell].
\]

For $s\in[0,s_{0}]$, we have
\[
\begin{aligned}
    \mathcal C_{2}'(s)
    ={}
    \rho e^{\frac{L_{3}}{2}s}
    \bigg[
        &\frac{L_{3}}{2}\cos(\beta s_{0}-\phi) - \beta\sin(\beta s_{0}-\phi) + \frac{L_{3}\beta}{2}(s_{0}-s)
        \sin(\beta s_{0}-\phi)
    \bigg].
\end{aligned}
\]
The expression in square brackets is largest at $s=0$. Since
\[
    s_{0}=\frac{2}{L_{3}}-\varepsilon,
\]
its value at $s=0$ is
\[
    \frac{L_{3}}{2}
    \left[
        \cos(\beta s_{0}-\phi) - \beta\varepsilon
        \sin(\beta s_{0}-\phi)
    \right],
\]
which is negative because
\[
    \tan(\beta s_{0}-\phi) > \frac{1}{\beta\varepsilon}.
\]
Therefore,
\[
    \mathcal C_{2}'(s)<0
    \qquad
    \text{for every }s\in[0,s_{0}].
\]
Since $\mathcal C_{2}(s_{0})>0$, it follows that
\[
    \mathcal C_{2}(s)>0
    \qquad
    \text{for every }s\in[0,s_{0}].
\]

On the first interval,
\[
    \mathcal C_{2}''
    -L_{3}\mathcal C_{2}'
    +\frac{L_{3}^{2}}{4}\mathcal C_{2}
    +L_{1}
    =0.
\]
Therefore, by \eqref{subsolution-discriminant},
\[
\begin{aligned}
    \mathcal C_{2}''
    -L_{3}\mathcal C_{2}'
    +L_{2}\mathcal C_{2}
    +L_{1}
    =
    -\left(
        \frac{L_{3}^{2}}{4}-L_{2}
    \right)\mathcal C_{2}\leq0.
\end{aligned}
\]
On the second interval, direct differentiation gives
\[
\begin{aligned}
    \mathcal C_2''-L_3\mathcal C_2'+L_2\mathcal C_2+L_1
    =-\left(\frac{L_3^2}{4}+\beta^2-L_2\right)\mathcal C_2
      -\frac{4L_1\beta^2}{L_3^2}<0.
\end{aligned}
\]
Thus \eqref{subsolution-control-ode} holds on each interval. Since $\mathcal C_2'$ is continuous at $s=s_0$, its distributional derivative has no singular part there. The inequality therefore holds weakly on $(0,\ell)$.

The second derivatives have a jump at the joining point, given by
\[
    \mathcal C_2''(s_0^+)-\mathcal C_2''(s_0^-)
    =-\rho\beta^2e^{L_3s_0/2}\cos(\beta s_0-\phi)<0.
\]
To obtain smooth global weights, we either use \eqref{critical-subsolution} or use \eqref{complete-subsolution-profile} with a cutoff supported away from the level $s=s_0$.

For either choice of $\mathcal C_2$, increase $\alpha_1$ in \eqref{C1} so that $\mathcal C_1(0)\geq\mathcal C_2(0)$. Since $\mathcal C_1$ is nondecreasing and $\mathcal C_2$ is nonincreasing, this gives $\mathcal C_1(s)\geq\mathcal C_2(s)$ throughout the interval. Consequently,
\begin{equation}\label{eq:local-weight-signs}
\begin{gathered}
    k_1\geq k_2>0,\\
    \partial_tk_1+\Delta_{\mu(t)}^{\nabla u}k_1-\mathcal Rk_1\geq0,
    \qquad
    \partial_tk_2+\Delta_{\mu(t)}^{\nabla u}k_2-\mathcal Rk_2\leq0
\end{gathered}
\end{equation}
on each regular region where both functions are smooth.

\subsection{Positive global weights}

Choose a nonzero smooth space-time cutoff $0\leq\chi\leq1$ supported in a region where $k_1,k_2$ are smooth and satisfy \eqref{eq:local-weight-signs}. Its support is compactly contained in this region relative to $M\times[t_0,b]$, so it may meet either endpoint in time. For each $t$, require
\begin{equation}\label{eq:original-cutoff-support}
    \operatorname{supp}_x\chi(\cdot,t)\Subset M_u(t)\cap M_p^{\mathrm{reg}}.
\end{equation}
Take $\chi=1$ on a smaller regular region. This construction assumes that $u$ has a nonempty smooth noncritical region. If $u$ is spatially constant throughout the interval, all spatial energy terms vanish and the radial construction is unnecessary.

Extend the local functions by setting
\begin{equation}\label{eq:independent-completions}
    \widetilde K_i=1-\chi+\chi k_i,\qquad i=1,2.
\end{equation}
The function $\chi(k_i-1)$ extends smoothly by zero outside the regular region. Since $\widetilde K_i$ is a convex combination of $1$ and $k_i$, we have $\widetilde K_1\geq\widetilde K_2>0$, with a uniform positive lower bound on $M\times[t_0,b]$. Both functions equal $1$ near the critical set.

Applying the product rule to $\widetilde K_i$ gives
\begin{equation}\label{eq:local-sign-completion-identity}
\begin{aligned}
    &\partial_t\widetilde K_i+\Delta_{\mu(t)}^{\nabla u}\widetilde K_i
           -\mathcal R\widetilde K_i\\
    &\quad=\chi\left(\partial_tk_i+\Delta_{\mu(t)}^{\nabla u}k_i-\mathcal Rk_i\right)\\
    &\qquad+(k_i-1)\left(\partial_t\chi+\Delta_{\mu(t)}^{\nabla u}\chi\right)
         +2d\chi(\nabla^{\nabla u}k_i)-\mathcal R(1-\chi).
\end{aligned}
\end{equation}
The first term has the sign given by \eqref{eq:local-weight-signs}. To control the terms arising from the cutoff and the constant density, choose $A\geq0$ such that
\begin{equation}\label{eq:cutoff-only-A}
\begin{aligned}
    A\geq\max_{i=1,2}\left\|
    \frac{(k_i-1)(\partial_t\chi+\Delta_{\mu(t)}^{\nabla u}\chi)
         +2d\chi(\nabla^{\nabla u}k_i)-\mathcal R(1-\chi)}
         {\widetilde K_i}
    \right\|_{L^\infty(M\times[t_0,b])}.
\end{aligned}
\end{equation}
The terms involving $k_i$ extend smoothly by zero because they are multiplied by derivatives of $\chi$ and supported in a compact noncritical region. The denominator is bounded below by a positive constant, and the numerator equals $-\mathcal R$ outside the cutoff support. Hence $A$ can be chosen finite. It depends on the local functions, the cutoff, and the coefficients of the equation, including the reference field $\nabla u$.

Set
\begin{equation}\label{eq:independent-global-weights}
    K_1=e^{A(t-t_0)}\widetilde K_1,
    \qquad
    K_2=e^{-A(t-t_0)}\widetilde K_2.
\end{equation}
Since $t\geq t_0$ and $\widetilde K_1\geq\widetilde K_2$, the weights are smooth, uniformly positive, and satisfy $K_1\geq K_2$. Combining \eqref{eq:local-weight-signs} with \eqref{eq:cutoff-only-A} gives
\begin{equation}\label{eq:global-super-from-local}
\begin{aligned}
    &e^{-A(t-t_0)}
      (\partial_tK_1+\Delta_{\mu(t)}^{\nabla u}K_1-\mathcal R K_1)\\
    &\quad=A\widetilde K_1
      +\chi(\partial_tk_1+\Delta_{\mu(t)}^{\nabla u}k_1-\mathcal Rk_1)\\
    &\qquad+(k_1-1)(\partial_t\chi+\Delta_{\mu(t)}^{\nabla u}\chi)
      +2d\chi(\nabla^{\nabla u}k_1)-\mathcal R(1-\chi)\geq0,
\end{aligned}
\end{equation}
and
\begin{equation}\label{eq:global-sub-from-local}
\begin{aligned}
    &e^{A(t-t_0)}
      (\partial_tK_2+\Delta_{\mu(t)}^{\nabla u}K_2-\mathcal R K_2)\\
    &\quad=-A\widetilde K_2
      +\chi(\partial_tk_2+\Delta_{\mu(t)}^{\nabla u}k_2-\mathcal Rk_2)\\
    &\qquad+(k_2-1)(\partial_t\chi+\Delta_{\mu(t)}^{\nabla u}\chi)
      +2d\chi(\nabla^{\nabla u}k_2)-\mathcal R(1-\chi)\leq0.
\end{aligned}
\end{equation}
Thus $K_1$ and $K_2$ are, respectively, a supersolution and a subsolution of the conjugate equation.

Extend $\nabla^{\nabla u}K_i$ by zero on $M\setminus M_u(t)$. The resulting vector field is smooth because $K_i$ is spatially constant near the critical set and its spatial derivatives are supported where $u$ is smooth and $du\ne0$. We define $\Delta_{\mu(t)}^{\nabla u}K_i$ as the divergence of this vector field. At a time when the cutoff is nonzero, choose a relatively compact open set $\Omega$ with smooth boundary such that $\operatorname{supp}_x\chi(\cdot,t)\Subset\Omega$ and $\overline\Omega$ lies in the regular region. For every smooth test function $\phi$,
\begin{equation}\label{eq:independent-global-ibp}
\begin{aligned}
    \int_M\phi\Delta_{\mu(t)}^{\nabla u}K_i\,d\mu(t)
    &=\int_\Omega\phi\Delta_{\mu(t)}^{\nabla u}K_i\,d\mu(t)=-\int_\Omega d\phi(\nabla^{\nabla u}K_i)\,d\mu(t)\\
     &=-\int_Md\phi(\nabla^{\nabla u}K_i)\,d\mu(t).
\end{aligned}
\end{equation}
The boundary term vanishes because the vector field is zero near $\partial\Omega$. If $\chi(\cdot,t)\equiv0$, every term in \eqref{eq:independent-global-ibp} is zero. Near the critical set, the left-hand sides of the conjugate inequalities reduce to $(A-\mathcal R)K_1$ and $(-A-\mathcal R)K_2$, which have the required signs by \eqref{eq:cutoff-only-A}. The inequalities therefore hold globally in the weak sense. Under the regularity assumptions stated below (see \eqref{eq:independent-regularity}), $u^2$ and $F^2(\nabla u)$ are admissible test functions for almost every time.

The ratio of the two weights is
\begin{equation}\label{eq:independent-weight-ratio}
    \frac{K_1}{K_2}=e^{2A(t-t_0)}
       \frac{1-\chi+\chi k_1}{1-\chi+\chi k_2}.
\end{equation}
This ratio may depend on space. We account for this dependence by a correction term $J$ in the frequency function (see \eqref{J-cf} in Section~\ref{sec:independent-frequency}). Set $d\mathcal V_i(t)=K_i\,d\mu(t)$. Since the exponential factors in \eqref{eq:independent-global-weights} are spatially constant, they do not affect the weighted Ricci curvature. For an affinely parametrized $F_t$-geodesic $\gamma$ with $\dot\gamma(0)=v\ne0$, the change-of-measure formula gives
\begin{equation}\label{eq:independent-curvature-change}
    \operatorname{Ric}_{K_i\mu(t)}^\infty(v)-\operatorname{Ric}_{\mu(t)}^\infty(v)
    =-\left.\frac{d^2}{ds^2}\log\widetilde K_i(\gamma(s),t)\right|_{s=0}.
\end{equation}
The integral bound on the metric variation and the change in weighted Ricci curvature is imposed on these global measures $d\mathcal V_i(t)$ (see \eqref{hv2mu<}) and is an additional assumption on the completed weights.

For the weighted energy identities used in the frequency estimates (see Section~\ref{sec:independent-frequency}), we assume
\begin{equation}\label{eq:independent-regularity}
\begin{gathered}
    u\in H^1((t_0,b);H^1(M))\cap L^2((t_0,b);H^2(M)),
    \qquad du\in C(M\times[t_0,b]),\\
    \Delta_{\mu(t)}u\in L^2((t_0,b);H^1(M)),
\end{gathered}
\end{equation}
together with the smooth noncritical neighborhood of the cutoff support chosen above. Sobolev norms may be defined using a fixed smooth auxiliary Riemannian metric. These assumptions imply $u\in C([t_0,b];H^1(M))$ and absolute continuity of the weighted energies
\[
    \int_M u^2 K_i\,d\mu(t),\qquad
    \int_M F_t^2(\nabla u)K_i\,d\mu(t).
\]
For almost every time, $u(\cdot,t)\in H^2(M)\cap C^1(M)$ and
\[
    \Delta_{\mathcal V_i(t)}u
    =\Delta_{\mu(t)}u+d\log K_i(\nabla u)\in H^1(M).
\]
Indeed, the additional term is smooth on the compact regular support of the spatial derivatives of $K_i$ and vanishes near the critical set. The integrated Bochner formula \cite[Theorem~3.6]{OhtaSturm2014Bochner} therefore applies with measure $d\mathcal V_i(t)$ for almost every time. Moreover, $u^2,F_t^2(\nabla u)\in H^1(M)\cap L^\infty(M)$ at those times, which justifies the tests used above.

\section{Parabolic frequency under coupled flows}\label{sec:independent-frequency}

Let $K_1,K_2$ be the global weights defined in \eqref{eq:independent-global-weights}, and set
\begin{equation}\label{dV1V2-cf}
    d\mathcal V_1=K_1\,d\mu(t),\qquad d\mathcal V_2=K_2\,d\mu(t),\qquad K_1\geq K_2>0.
\end{equation}
These are smooth positive measures on $M$, and their ratio may depend on both space and time. We first consider the heat equation
\begin{equation}\label{H3}
    u_t-\Delta_{\mu(t)}u=0.
\end{equation}
Throughout this section, $u$ satisfies the regularity assumptions of Section~\ref{sec:local-conjugate-construction}, which justify differentiation of the energies and application of the integrated Bochner formula. Define
\begin{equation}\label{ID-cf}
    I=\int_Mu^2\,d\mathcal V_1,\qquad
    D=-\tau\int_MF^2(\nabla u)\,d\mathcal V_2,\qquad \tau=t_1-t>0.
\end{equation}
On intervals where $I>0$, let
\begin{equation}\label{J-cf}
    J(t)=\max_M\left[\left(1-\frac{K_1}{K_2}\right)\frac D{\tau I}\right]
       =\left(\max_M\frac{K_1}{K_2}-1\right)\frac{-D}{\tau I}\geq0,
\end{equation}
and
\begin{equation}\label{UwithJ}
    U(t)=\frac{D(t)}{I(t)}
      \exp\left(\int_a^t\left[\frac{1-\kappa(s)}{\tau(s)}-2J(s)\right]ds\right).
\end{equation}
Here $a$ belongs to the interval on which $I>0$, and $\kappa$ is a bounded function controlling the integral of the metric variation and the change in weighted Ricci curvature (see \eqref{hv2mu<}). One may take $a=t_0$ if $I(t_0)>0$.

\begin{theorem}\label{thm:independent-heat}
Let $(M,F(t),\mu(t))$, $t\in[t_0,b]$ with $b<t_1$, be a smooth family of closed Finsler metric measure spaces satisfying \eqref{F2}. Let $u$ solve \eqref{H3} with the regularity specified above, and let $K_1,K_2$ be the weights constructed in Section~\ref{sec:local-conjugate-construction}. Suppose that a bounded function $\kappa$ satisfies
\begin{equation}\label{hv2mu<}
\begin{aligned}
    &\int_M\left[h(\nabla u)+\operatorname{Ric}_{\mathcal V_2}^\infty(\nabla u)
        -\operatorname{Ric}_{\mu(t)}^\infty(\nabla u)\right]d\mathcal V_2\\
    &\qquad\leq\frac{\kappa(t)}{2\tau(t)}
           \int_MF^2(\nabla u)\,d\mathcal V_2
       =-\frac{\kappa(t)}{2\tau(t)^2}D(t)
\end{aligned}
\end{equation}
for almost every $t\in[t_0,b]$. Then the frequency \eqref{UwithJ} is nondecreasing on every interval on which $I>0$.
\end{theorem}

\begin{proof}
The volume-form derivatives are
\[
    \partial_t(d\mathcal V_i)
      =(\partial_tK_i-\mathcal R K_i)\,d\mu(t).
\]
The supersolution and subsolution inequalities give
\[
    \partial_tK_1-\mathcal R K_1\geq-\Delta_{\mu(t)}^{\nabla u}K_1,
    \qquad
    \partial_tK_2-\mathcal R K_2\leq-\Delta_{\mu(t)}^{\nabla u}K_2.
\]
We use these inequalities in their global weak form, as established in Section~\ref{sec:local-conjugate-construction}.

Differentiating $I$ and using the supersolution inequality, we obtain
\begin{equation}\label{eq:independent-I-heat}
\begin{aligned}
    I'={}&2\int_Mu u_t\,d\mathcal V_1
       +\int_Mu^2(\partial_tK_1-\mathcal R K_1)\,d\mu(t)\\
    \geq{}&2\int_Mu\Delta_{\mu(t)}u\,K_1\,d\mu(t)
       -\int_Mu^2\Delta_{\mu(t)}^{\nabla u}K_1\,d\mu(t)\\
    ={}&-2\int_Md(uK_1)(\nabla u)\,d\mu(t)
       +2\int_Mu\,dK_1(\nabla u)\,d\mu(t)\\
    ={}&-2\int_MF^2(\nabla u)\,d\mathcal V_1.
\end{aligned}
\end{equation}
Since $\tau'=-1$, differentiation of $D$ gives
\begin{equation}\label{D'=}
\begin{aligned}
    D'={}&\int_MF^2(\nabla u)\,d\mathcal V_2
        -\tau\int_M\partial_t[F^2(\nabla u)]\,d\mathcal V_2\\
       &-\tau\int_MF^2(\nabla u)(\partial_tK_2-\mathcal R K_2)\,d\mu(t)\\
    ={}&-\frac D\tau-\tau\int_M\partial_t[F^2(\nabla u)]\,d\mathcal V_2
       -\tau\int_MF^2(\nabla u)(\partial_tK_2-\mathcal R K_2)\,d\mu(t).
\end{aligned}
\end{equation}
Lemma~\ref{lem1} and the subsolution inequality therefore imply
\begin{equation}\label{ineq-D'+dtau}
\begin{aligned}
    D'+\frac D\tau
    ={}&-2\tau\int_M[du_t(\nabla u)+h(\nabla u)]\,d\mathcal V_2\\
       &-\tau\int_MF^2(\nabla u)(\partial_tK_2-\mathcal R K_2)\,d\mu(t)\\
    \geq{}&-2\tau\int_M[du_t(\nabla u)+h(\nabla u)]\,d\mathcal V_2
       +\tau\int_MF^2(\nabla u)\Delta_{\mu(t)}^{\nabla u}K_2\,d\mu(t)\\
    ={}&-2\tau\int_M[du_t(\nabla u)+h(\nabla u)]\,d\mathcal V_2
       -\tau\int_Md[F^2(\nabla u)](\nabla^{\nabla u}K_2)\,d\mu(t).
\end{aligned}
\end{equation}
The last equality follows from the weak divergence identity with test function $F^2(\nabla u)$. This is justified by the stated regularity and the compact support of $\nabla^{\nabla u}K_2$ in the smooth noncritical region.

On $M_u$, the Bochner--Weitzenb\"ock formula gives
\[
\begin{aligned}
    \partial_tF^2(\nabla u)-\Delta_{\mu(t)}^{\nabla u}F^2(\nabla u)
      ={}&2du_t(\nabla u)+2h(\nabla u)-2d(\Delta_{\mu(t)}u)(\nabla u)\\
       &-2\operatorname{Ric}_{\mu(t)}^\infty(\nabla u)
        -2\|\nabla^2u\|_{HS(\nabla u)}^2\\
      ={}&2\left[h(\nabla u)-\operatorname{Ric}_{\mu(t)}^\infty(\nabla u)
        -\|\nabla^2u\|_{HS(\nabla u)}^2\right].
\end{aligned}
\]
Applying its integrated form \cite[Theorem~3.6]{OhtaSturm2014Bochner} with test function $K_2$, we have
\[
\begin{aligned}
    -\frac12\int_Md[F^2(\nabla u)](\nabla^{\nabla u}K_2)\,d\mu(t)
      =\int_M\left[d(\Delta_{\mu(t)}u)(\nabla u)
       +\operatorname{Ric}_{\mu(t)}^\infty(\nabla u)
       +\|\nabla^2u\|_{HS(\nabla u)}^2\right]d\mathcal V_2.
\end{aligned}
\]
Since $u_t=\Delta_{\mu(t)}u$, substitution into \eqref{ineq-D'+dtau} yields
\begin{equation}\label{ineq-D'+Dtau-2}
    D'+\frac D\tau\geq
    2\tau\int_M\left[\operatorname{Ric}_{\mu(t)}^\infty(\nabla u)
      +\|\nabla^2u\|_{HS(\nabla u)}^2-h(\nabla u)\right]d\mathcal V_2.
\end{equation}
Applying the integrated Bochner formula with measure $d\mathcal V_2$ and test function $1$, followed by integration by parts, gives
\begin{equation}\label{eq-intnabla2u}
\begin{aligned}
    \int_M\|\nabla^2u\|_{HS(\nabla u)}^2\,d\mathcal V_2
    &=-\int_Md(\Delta_{\mathcal V_2}u)(\nabla u)\,d\mathcal V_2
       -\int_M\operatorname{Ric}_{\mathcal V_2}^\infty(\nabla u)\,d\mathcal V_2\\
    &=\int_M\left[(\Delta_{\mathcal V_2}u)^2
       -\operatorname{Ric}_{\mathcal V_2}^\infty(\nabla u)\right]d\mathcal V_2.
\end{aligned}
\end{equation}
The Hessian is independent of the measure. Substituting \eqref{eq-intnabla2u} into \eqref{ineq-D'+Dtau-2} and applying \eqref{hv2mu<}, we obtain
\begin{equation}\label{eq:independent-curvature-estimate}
\begin{aligned}
    D'+\frac D\tau
    \geq{}&2\tau\int_M(\Delta_{\mathcal V_2}u)^2\,d\mathcal V_2-2\tau\int_M\left[h(\nabla u)+\operatorname{Ric}_{\mathcal V_2}^\infty(\nabla u)
         -\operatorname{Ric}_{\mu(t)}^\infty(\nabla u)\right]d\mathcal V_2\\
    \geq{}&2\tau\int_M(\Delta_{\mathcal V_2}u)^2\,d\mathcal V_2
       -\kappa\int_MF^2(\nabla u)\,d\mathcal V_2\\
    ={}&2\tau\int_M(\Delta_{\mathcal V_2}u)^2\,d\mathcal V_2+\frac\kappa\tau D.
\end{aligned}
\end{equation}

To compare the energies computed with the two measures, we use the definition of $J$, which gives
\[
    JI=\left(\max_M\frac{K_1}{K_2}-1\right)
        \int_MF^2(\nabla u)\,d\mathcal V_2,
\]
and
\[
    \int_MF^2(\nabla u)\,d\mathcal V_1
    \leq\max_M\frac{K_1}{K_2}\int_MF^2(\nabla u)\,d\mathcal V_2.
\]
Since $D\leq0$, these imply
\begin{equation}\label{-2JID}
    -2JDI+2D\int_MF^2(\nabla u)\,d\mathcal V_1\geq-\frac{2D^2}{\tau}.
\end{equation}
Differentiating \eqref{UwithJ} and using \eqref{eq:independent-curvature-estimate} and \eqref{-2JID}, we obtain
\begin{equation}\label{eq:independent-heat-quotient}
\begin{aligned}
    &e^{-\int_a^t[(1-\kappa(s))/\tau(s)-2J(s)]\,ds}I^2U'\\
    &\quad=ID'-DI'+\left(\frac{1-\kappa}{\tau}-2J\right)DI\\
    &\quad\geq2\tau I\int_M(\Delta_{\mathcal V_2}u)^2\,d\mathcal V_2
         +2D\int_MF^2(\nabla u)\,d\mathcal V_1-2JDI\\
    &\quad\geq2\tau I\int_M(\Delta_{\mathcal V_2}u)^2\,d\mathcal V_2
         -\frac{2D^2}{\tau}\geq0.
\end{aligned}
\end{equation}
The final step follows from integration by parts and the Cauchy--Schwarz inequality, since
\begin{equation}\label{eq:independent-cauchy-schwarz}
\begin{aligned}
    \frac{D^2}{\tau^2}
    &=\left(\int_Mu\Delta_{\mathcal V_2}u\,d\mathcal V_2\right)^2\\
    &\leq\left(\int_Mu^2\,d\mathcal V_2\right)
          \left(\int_M(\Delta_{\mathcal V_2}u)^2\,d\mathcal V_2\right)\\
    &\leq I\int_M(\Delta_{\mathcal V_2}u)^2\,d\mathcal V_2.
\end{aligned}
\end{equation}
The last inequality follows from $K_1\geq K_2$, completing the proof.
\end{proof}

A nonzero spatially constant heat solution satisfies $D=U=0$ for any choice of weights, so equality $U'=0$ alone does not imply $K_1=K_2$. If $K_1=K_2$ is an exact solution of the conjugate equation, then $J=0$, and equality holds precisely when equality holds in \eqref{hv2mu<} and $\Delta_{\mathcal V_2}u=c(t)u$, where
\[
    c(t)=\frac{D(t)}{\tau(t)I(t)}
       =\frac{U(t)}{\tau(t)}
          \exp\left(-\int_a^t\frac{1-\kappa(s)}{\tau(s)}\,ds\right).
\]

We next consider
\begin{equation}\label{H4}
    u_t-\Delta_{\mu(t)}u=Q(u,\nabla u),\qquad
    |Q(u,\nabla u)|\leq C(t)(|u|+F(\nabla u)),
\end{equation}
where $C$ is bounded and nonnegative. We obtain a differential inequality for the frequency and a backward uniqueness result.

\begin{theorem}\label{thm:independent-backward}
Let $(M,F(t),\mu(t))$ be a smooth family of closed Finsler metric measure spaces satisfying \eqref{F2} on $[t_0,b]$, where $b<t_1$. Let $u$ solve \eqref{H4} with the regularity specified above, and let $K_1,K_2$ be the weights constructed in Section~\ref{sec:local-conjugate-construction}. Suppose that a bounded function $\kappa$ satisfies \eqref{hv2mu<}. Then, on every interval on which $I>0$, the frequency satisfies
\begin{equation}\label{U'geq}
\begin{aligned}
    U'\geq{}&C(t)\left[C(t)+(3+J(t))\max_M\frac{K_1}{K_2}-3\right]U\\
       &-\tau(t)C(t)^2
          \exp\left(\int_a^t\left[\frac{1-\kappa(s)}{\tau(s)}-2J(s)\right]ds\right).
\end{aligned}
\end{equation}
Moreover, if $u(\cdot,b)=0$, then $u\equiv0$ on $M\times[t_0,b]$.
\end{theorem}

\begin{proof}
Fix the weights $K_1,K_2$ on $M\times[t_0,b]$. As in the preceding proof, the supersolution inequality gives
\begin{equation}\label{I'-cf}
\begin{aligned}
    I'={}&2\int_Mu u_t\,d\mathcal V_1
       +\int_Mu^2(\partial_tK_1-\mathcal R K_1)\,d\mu(t)\\
    \geq{}&2\int_Mu(\Delta_{\mu(t)}u+Q)K_1\,d\mu(t)
       -\int_Mu^2\Delta_{\mu(t)}^{\nabla u}K_1\,d\mu(t)\\
    ={}&2\int_M[uQ-F^2(\nabla u)]\,d\mathcal V_1.
\end{aligned}
\end{equation}
Multiplication by $-D\geq0$ gives
\begin{equation}\label{-DI'}
    -DI'\geq2D\int_MF^2(\nabla u)\,d\mathcal V_1
        -2D\int_MuQ\,d\mathcal V_1.
\end{equation}
Combining this with \eqref{-2JID} yields
\begin{equation}\label{-2JID-DI'}
    -2JDI-DI'\geq-\frac{2D^2}{\tau}-2D\int_MuQ\,d\mathcal V_1.
\end{equation}
For the numerator, the global subsolution inequality and the integrated Bochner formula give
\[
\begin{aligned}
    D'+\frac D\tau\geq{}&-2\tau\int_M[du_t(\nabla u)+h(\nabla u)]\,d\mathcal V_2
        -\tau\int_Md[F^2(\nabla u)](\nabla^{\nabla u}K_2)\,d\mu(t)\\
    ={}&2\tau\int_M\left[-dQ(\nabla u)-h(\nabla u)
        +\operatorname{Ric}_{\mu(t)}^\infty(\nabla u)
        +\|\nabla^2u\|_{HS(\nabla u)}^2\right]d\mathcal V_2.
\end{aligned}
\]
The regularity assumptions give $Q=u_t-\Delta_{\mu(t)}u\in H^1(M)$ for almost every $t$, so integration by parts yields
\[
    -\int_MdQ(\nabla u)\,d\mathcal V_2
       =\int_MQ\Delta_{\mathcal V_2}u\,d\mathcal V_2.
\]
Using \eqref{eq-intnabla2u} and \eqref{hv2mu<}, we conclude
\begin{equation}\label{D'>-cf}
\begin{aligned}
    D'+\frac D\tau
    \geq{}&2\tau\int_M\left[(\Delta_{\mathcal V_2}u)^2
          +Q\Delta_{\mathcal V_2}u\right]d\mathcal V_2\\
       &-2\tau\int_M\left[h(\nabla u)+\operatorname{Ric}_{\mathcal V_2}^\infty(\nabla u)
          -\operatorname{Ric}_{\mu(t)}^\infty(\nabla u)\right]d\mathcal V_2\\
    \geq{}&2\tau\int_M\left[(\Delta_{\mathcal V_2}u)^2
          +Q\Delta_{\mathcal V_2}u\right]d\mathcal V_2+\frac\kappa\tau D.
\end{aligned}
\end{equation}
The quotient derivative consequently satisfies
\begin{equation}\label{eq:independent-source-quotient}
\begin{aligned}
    &e^{-\int_a^t[(1-\kappa(s))/\tau(s)-2J(s)]\,ds}I^2U'\\
    &\quad\geq2\tau\Bigg[I\int_M\left[(\Delta_{\mathcal V_2}u)^2
         +Q\Delta_{\mathcal V_2}u\right]d\mathcal V_2
         -\frac{D^2}{\tau^2}-\frac D\tau\int_MuQ\,d\mathcal V_1\Bigg].
\end{aligned}
\end{equation}
Since $D/\tau=\int_Mu\Delta_{\mathcal V_2}u\,d\mathcal V_2$, completing the square in the bracketed expression gives
\begin{align*}
    &I\int_M\left[(\Delta_{\mathcal V_2}u)^2+Q\Delta_{\mathcal V_2}u\right]d\mathcal V_2
         -\frac{D^2}{\tau^2}-\frac D\tau\int_MuQ\,d\mathcal V_1\\
    &=I\int_M\left(\Delta_{\mathcal V_2}u+\frac12Q\right)^2d\mathcal V_2
       -\left[\int_Mu\left(\Delta_{\mathcal V_2}u+\frac12Q\right)d\mathcal V_2\right]^2\\
    &\quad-\frac I4\int_MQ^2\,d\mathcal V_2
       +\frac14\left(\int_MuQ\,d\mathcal V_2\right)^2\\
    &\quad+\int_MF^2(\nabla u)\,d\mathcal V_2
          \int_MuQ\left(\frac{K_1}{K_2}-1\right)d\mathcal V_2\\
    &\geq-\frac I4\int_MQ^2\,d\mathcal V_2
       +\int_MF^2(\nabla u)\,d\mathcal V_2
          \int_MuQ\left(\frac{K_1}{K_2}-1\right)d\mathcal V_2.
\end{align*}
The first two terms after the equality have nonnegative sum by Cauchy--Schwarz and $\int_Mu^2\,d\mathcal V_2\leq I$, while $\frac14(\int_MuQ\,d\mathcal V_2)^2$ is also nonnegative.

The growth bound on $Q$ gives
\[
    \frac12\int_MQ^2\,d\mathcal V_2
       \leq C(t)^2\int_M[u^2+F^2(\nabla u)]\,d\mathcal V_2
       \leq C(t)^2\left(I-\frac D\tau\right),
\]
and
\[
    \int_MuQ\left(\frac{K_1}{K_2}-1\right)d\mathcal V_2
       \geq-\frac{C(t)}2\int_M
          \left(\frac{K_1}{K_2}-1\right)[3u^2+F^2(\nabla u)]\,d\mathcal V_2.
\]
For the gradient part,
\[
    \frac1I\int_MF^2(\nabla u)
       \left(\frac{K_1}{K_2}-1\right)d\mathcal V_2\leq J(t),
\]
while the corresponding $u^2$ integral divided by $I$ is at most $\max_M(K_1/K_2)-1$. Consequently,
\[
\begin{aligned}
    U'\geq{}&\left[C(t)^2+3C(t)\left(\max_M\frac{K_1}{K_2}-1\right)+C(t)J(t)\right]U\\
       &-\tau C(t)^2e^{\int_a^t[(1-\kappa(s))/\tau(s)-2J(s)]\,ds}.
\end{aligned}
\]
Since $U\leq0$, $J\geq0$, and $\max_M(K_1/K_2)\geq1$, this implies \eqref{U'geq}.

To prove backward uniqueness, we first estimate the energy quotient with both integrals taken against $d\mathcal V_2$. This will also bound $J$ near a possible zero of $I$.

Since $K_2$ is a smooth positive subsolution, its normalized residual is nonpositive and bounded on $M\times[t_0,b]$. Choose $\Lambda\geq1$ so that, for almost every $t$,
\begin{equation}\label{Lambda-cf}
\begin{gathered}
    \frac{|\kappa(t)|}{\tau(t)}+C(t)^2+
    \left\|\frac{\partial_tK_2+\Delta_{\mu(t)}^{\nabla u}K_2-\mathcal R K_2}{K_2}
         \right\|_{L^\infty(M)}\leq\Lambda,\\
    C(t)\leq\Lambda,\qquad \max_M\frac{K_1}{K_2}\leq\Lambda.
\end{gathered}
\end{equation}
Such a constant exists because $b<t_1$, $C$ and $\kappa$ are bounded, and the smooth weights have positive lower bounds, with their spatial derivatives supported in the smooth noncritical region.

Differentiation and integration by parts give
\begin{equation}\label{eq:second-measure-mass-identity}
\begin{aligned}
    \frac{d}{dt}\int_Mu^2\,d\mathcal V_2
    ={}&\frac{2D}{\tau}+2\int_MuQ\,d\mathcal V_2\\
       &+\int_Mu^2
       \frac{\partial_tK_2+\Delta_{\mu(t)}^{\nabla u}K_2-\mathcal R K_2}{K_2}
          \,d\mathcal V_2.
\end{aligned}
\end{equation}
Combining this identity with \eqref{D'>-cf}, differentiating the quotient, and using $D/\tau=\int_Mu\Delta_{\mathcal V_2}u\,d\mathcal V_2$, we obtain
\begin{align}
    &\left(\frac{D}{\tau\int_Mu^2\,d\mathcal V_2}\right)'
       \notag\\
    &\quad\geq\left(\frac\kappa\tau+
       \left\|\frac{\partial_tK_2+\Delta_{\mu(t)}^{\nabla u}K_2-\mathcal R K_2}{K_2}
          \right\|_{L^\infty(M)}\right)
          \frac{D}{\tau\int_Mu^2\,d\mathcal V_2}
       \notag\\
    &\qquad+\frac{2}{\int_Mu^2\,d\mathcal V_2}
       \int_M\left[(\Delta_{\mathcal V_2}u)^2+Q\Delta_{\mathcal V_2}u\right]d\mathcal V_2
       \notag\\
    &\qquad-\frac{2}{(\int_Mu^2\,d\mathcal V_2)^2}
       \left[\frac{D^2}{\tau^2}+\frac D\tau\int_MuQ\,d\mathcal V_2\right]
       \notag\\
    &\quad\geq\left(\frac\kappa\tau+
       \left\|\frac{\partial_tK_2+\Delta_{\mu(t)}^{\nabla u}K_2-\mathcal R K_2}{K_2}
          \right\|_{L^\infty(M)}\right)
          \frac{D}{\tau\int_Mu^2\,d\mathcal V_2}
       -\frac{\int_MQ^2\,d\mathcal V_2}{2\int_Mu^2\,d\mathcal V_2}
       \notag\\
    &\quad\geq\left(\frac\kappa\tau+C(t)^2+
       \left\|\frac{\partial_tK_2+\Delta_{\mu(t)}^{\nabla u}K_2-\mathcal R K_2}{K_2}
          \right\|_{L^\infty(M)}\right)
          \frac{D}{\tau\int_Mu^2\,d\mathcal V_2}-C(t)^2
       \notag\\
    &\quad\geq\Lambda\frac{D}{\tau\int_Mu^2\,d\mathcal V_2}-\Lambda^2.
       \label{eq:single-measure-quotient-estimate}
\end{align}
The second inequality follows by completing the square with $\Delta_{\mathcal V_2}u+Q/2$ and applying Cauchy--Schwarz in $L^2(M,d\mathcal V_2)$. The last inequality uses $D\leq0$ and \eqref{Lambda-cf}.

Suppose $u(\cdot,b)=0$ but $I(a)>0$ for some $a\in[t_0,b)$. Let $b_0$ be the first zero of $I$ after $a$. Both $I$ and $\int_Mu^2\,d\mathcal V_2$ are positive on $[a,b_0)$ because both densities are positive. Multiplying \eqref{eq:single-measure-quotient-estimate} by $e^{-\Lambda(t-a)}$ and integrating yields
\[
\begin{aligned}
    e^{-\Lambda(t-a)}\frac{D(t)}{\tau(t)\int_Mu(x,t)^2\,d\mathcal V_2(t)}
       -\frac{D(a)}{\tau(a)\int_Mu(x,a)^2\,d\mathcal V_2(a)}
    \geq-\Lambda\left(1-e^{-\Lambda(t-a)}\right).
\end{aligned}
\]
Consequently,
\begin{equation}\label{U-lower-Lambda}
\begin{aligned}
    \frac{D(t)}{\tau(t)\int_Mu(x,t)^2\,d\mathcal V_2(t)}
    \geq{}&e^{\Lambda(t-a)}
       \frac{D(a)}{\tau(a)\int_Mu(x,a)^2\,d\mathcal V_2(a)}\\
       &-\Lambda\left(e^{\Lambda(t-a)}-1\right).
\end{aligned}
\end{equation}
The right-hand side has a finite lower bound as $t\rightarrow b_0$. Since $I\geq\int_Mu^2\,d\mathcal V_2$, we also obtain
\begin{equation}\label{eq:J-endpoint-bound-proved}
\begin{aligned}
    0\leq J(t)
    &\leq(\Lambda-1)\frac{-D(t)}{\tau(t)\int_Mu(x,t)^2\,d\mathcal V_2(t)}\\
    &\leq(\Lambda-1)\left[
      e^{\Lambda(t-a)}\frac{-D(a)}{\tau(a)\int_Mu(x,a)^2\,d\mathcal V_2(a)}
       +\Lambda(e^{\Lambda(t-a)}-1)\right].
\end{aligned}
\end{equation}
Thus $J$ is bounded on $[a,b_0)$ and hence integrable up to $b_0$.

Finally, \eqref{I'-cf} and Young's inequality give
\begin{equation}\label{logI'geq}
\begin{aligned}
    (\log I)'
    &\geq-3C(t)-\frac{2+C(t)}I\int_MF^2(\nabla u)\,d\mathcal V_1\\
    &\geq-3C(t)-(2+C(t))\max_M\frac{K_1}{K_2}
       \frac{\int_MF^2(\nabla u)\,d\mathcal V_2}{\int_Mu^2\,d\mathcal V_2}\\
    &\geq(2+\Lambda)\Lambda
       \frac{D}{\tau\int_Mu^2\,d\mathcal V_2}-3\Lambda.
\end{aligned}
\end{equation}
Substituting \eqref{U-lower-Lambda} and integrating from $a$ to $s<b_0$ gives
\begin{equation}\label{Ia<Ib}
\begin{aligned}
    I(s)\geq I(a)\exp\Bigg\{&
      (2+\Lambda)\left[
       \frac{D(a)}{\tau(a)\int_Mu(x,a)^2\,d\mathcal V_2(a)}-\Lambda\right]
           (e^{\Lambda(s-a)}-1)\\
      &+\left[(2+\Lambda)\Lambda^2-3\Lambda\right](s-a)\Bigg\}>0.
\end{aligned}
\end{equation}
The right-hand side tends to a strictly positive number as $s\rightarrow b_0$, contradicting $I(b_0)=0$. Therefore $u\equiv0$ on $M\times[t_0,b]$. If $u$ is spatially constant throughout the interval, the equation reduces to $u_t=Q(u,0)$ with $|u_t|\leq C(t)|u|$, and the conclusion follows directly from Gronwall's inequality.
\end{proof}

\clearpage
\appendix
\section{Variational existence along a coupled flow}\label{app:heat-existence}

We give a detailed proof of Proposition~\ref{prop:coupled-heat-existence}. The argument uses the scalar $L^2$ structure and the convex Finsler energy, as in the static construction of Ohta--Sturm \cite[Section~3]{OhtaSturm2009HeatFlow}. We first state the time-dependent existence theorem of Kopfer \cite[Theorem~5.5]{Kopfer2018GradientFlow}, verify its hypotheses, and then describe the implicit time approximation and its convergence in the present compact setting.

\subsection{Time-dependent convex energies}\label{app:abstract-heat-flow}

Let $\mathcal H$ be a separable real Hilbert space with a family of inner products $\langle\cdot,\cdot\rangle_t$, $t\in[t_0,b]$, and write $\|z\|_t^2=\langle z,z\rangle_t$. For a convex functional $E_t:\mathcal H\to[0,+\infty]$, its subdifferential with respect to this inner product is
\begin{equation}\label{eq:app-subdifferential}
    \partial^{\,t}E_t(z)
    =
    \left\{\zeta\in\mathcal H\ \middle|\
       E_t(w)-E_t(z)\geq\langle\zeta,w-z\rangle_t
       \text{ for all }w\in\mathcal H\right\},
\end{equation}
where $E_t(z)<\infty$. The superscript $t$ specifies the inner product used in the subdifferential. A curve $z$ belongs to $AC^2([t_0,b];\mathcal H)$ if there is $m\in L^2(t_0,b)$ such that
\[
    \|z(t)-z(s)\|_{t_0}\leq\int_s^t m(r)\,dr
    \qquad(t_0\leq s\leq t\leq b).
\]
In a Hilbert space, these curves are precisely the continuous representatives of $H^1((t_0,b);\mathcal H)$, whose weak time derivatives lie in $L^2((t_0,b);\mathcal H)$.

\begin{theorem}[Kopfer {\cite[Theorem~5.5]{Kopfer2018GradientFlow}}]\label{thm:app-kopfer}
Suppose that the inner products satisfy
\begin{equation}\label{eq:app-hilbert-comparison}
    e^{-L|t-s|}\|z\|_s\leq\|z\|_t
       \leq e^{L|t-s|}\|z\|_s
    \qquad(z\in\mathcal H,\ s,t\in[t_0,b])
\end{equation}
for some $L\geq0$. Assume that each $E_t:\mathcal H\to[0,+\infty]$ is convex and lower semicontinuous, that its effective domain
\[
    D=\{z\in\mathcal H\mid E_t(z)<\infty\}
\]
is nonempty and independent of $t$, and that
\begin{equation}\label{eq:app-relative-time}
    |E_t(z)-E_s(z)|\leq C_1|t-s|E_t(z)
    \qquad(z\in D,\ s,t\in[t_0,b])
\end{equation}
for some $C_1\geq0$. Then every $z_0\in D$ determines a unique curve $z\in AC^2([t_0,b];\mathcal H)$ satisfying
\[
    z(t_0)=z_0,\qquad
    -z'(t)\in\partial^{\,t}E_t(z(t))
    \quad\text{for almost every }t\in(t_0,b).
\]
\end{theorem}

The assumptions are those at the beginning of Sections~5 and~5.1 of \cite{Kopfer2018GradientFlow}. The estimates in \cite[Proposition~5.3]{Kopfer2018GradientFlow} also give bounded energy along the constructed curve. The functional may be nonquadratic, and the initial datum is required to have finite energy, without an assumption that its subdifferential is nonempty.

\subsection{The Finsler energy and the heat equation}\label{app:finsler-heat-proof}

\begin{proof}[Proof of Proposition~\ref{prop:coupled-heat-existence}]
We work on the prescribed smooth family $(M,F(t),\mu(t))$, $t\in[t_0,b]$, and write $F_t=F(t)$. Fix a smooth auxiliary Riemannian metric $\widehat g$ for the Sobolev norms. Set
\begin{equation}\label{eq:app-l2-energy}
\begin{gathered}
    \mathcal H=L^2(M,\mu(t_0)),\qquad
    \langle f,k\rangle_t=\int_Mfk\,d\mu(t),\\
    E_t(f)=
    \begin{cases}
       \displaystyle\frac12\int_MF_t^{*2}(df)\,d\mu(t),
          &f\in H^1(M),\\[1ex]
       +\infty,&f\notin H^1(M).
    \end{cases}
\end{gathered}
\end{equation}
The inner product is the usual one on real-valued scalar functions. Thus $\mathcal H$ is Hilbert even when the tangent norms of $F_t$ do not come from inner products. The Finsler structure enters through the generally nonquadratic functional $E_t$.

\smallskip
\noindent\emph{The Hilbert norms and the common domain.}
Choose $A,B\geq0$ so that
\[
    |H(x,v,t)|\leq A,\qquad
    |\mathcal R(x,t)|\leq B
\]
on the compact time interval. The first bound is obtained on the unit tangent bundle of $\widehat g$ and extends to every $v\ne0$ because $H$ is zero homogeneous. If $d\mu(t)=\rho_t\,d\mu(t_0)$, the second equation in \eqref{F2} gives
\[
    \rho_t(x)=\exp\left(-\int_{t_0}^t\mathcal R(x,r)\,dr\right),
    \qquad
    e^{-B|t-s|}\,d\mu(s)\leq d\mu(t)
       \leq e^{B|t-s|}\,d\mu(s).
\]
In particular, all the measures have the same null sets and define the same $L^2$ vector space. Their norms satisfy
\begin{equation}\label{eq:app-l2-time}
    e^{-B|t-s|/2}\|f\|_s\leq\|f\|_t
       \leq e^{B|t-s|/2}\|f\|_s.
\end{equation}
The space $\mathcal H$ is separable, and the equivalent norms are complete. Thus \eqref{eq:app-hilbert-comparison} holds with $L=B/2$.

Positivity and continuity of $F_t^*$ on the unit cotangent bundle of $\widehat g$, together with homogeneity and compactness, give uniform equivalence with $|\cdot|_{\widehat g^{-1}}$. Combining this with the measure comparison, there is $c>0$ such that
\begin{equation}\label{eq:app-energy-comparison}
    c\int_M|df|_{\widehat g^{-1}}^2\,d\mu(t_0)
    \leq 2E_t(f)
    \leq c^{-1}\int_M|df|_{\widehat g^{-1}}^2\,d\mu(t_0).
\end{equation}
Hence the finite-energy domain in $\mathcal H$ is the fixed space $H^1(M)$. The energy controls the differential, while the ambient $L^2$ norm controls the function itself.

\smallskip
\noindent\emph{Convexity and lower semicontinuity.}
Write
\[
    e_t(x,\xi)=\frac12F_t^{*2}(x,\xi),\qquad
    J_t^*(x,\xi)=D_\xi e_t(x,\xi).
\]
The function $e_t(x,\cdot)$ is convex. Its derivative is the dual Legendre map, and $J_t^*(du)=\nabla u$ \cite[Lemma~1.1]{OhtaSturm2009HeatFlow}. Moreover,
\[
    e_t(x,\xi)\leq C|\xi|_{\widehat g^{-1}}^2,
    \qquad
    |J_t^*(x,\xi)|_{\widehat g}
       \leq C|\xi|_{\widehat g^{-1}},
\]
uniformly on the compact space-time region. The first estimate shows that $e_t$ is differentiable at $\xi=0$ with derivative zero. The second, obtained by homogeneity from the unit cotangent bundle, shows that this derivative is continuous there. Thus the first derivative of $e_t$ is defined even at zero covectors.

Since the differential is linear, pointwise convexity yields
\[
    E_t((1-\theta)f+\theta k)
       \leq(1-\theta)E_t(f)+\theta E_t(k),
       \qquad 0\leq\theta\leq1.
\]
Also $E_t\geq0$ and $E_t(0)=0$. To check lower semicontinuity in $\mathcal H$, suppose that $f_j\to f$ in $\mathcal H$ and that $\liminf_jE_t(f_j)<\infty$. Choose a subsequence attaining this lower limit. By \eqref{eq:app-energy-comparison} it is bounded in $H^1(M)$, so a further subsequence converges weakly to $f$ in $H^1(M)$. The supporting inequality for $e_t$ gives
\[
    E_t(f_j)\geq E_t(f)
       +\int_M(df_j-df)\bigl(J_t^*(df)\bigr)\,d\mu(t).
\]
The last integral tends to zero because $df_j\rightharpoonup df$ in $L^2$ and $J_t^*(df)\in L^2$. Therefore $E_t(f)\leq\liminf_jE_t(f_j)$. The same growth bound and the integral mean value formula give
\begin{equation}\label{eq:app-h1-continuity}
    |E_t(f)-E_t(k)|
       \leq C\bigl(\|df\|_{L^2}+\|dk\|_{L^2}\bigr)
                    \|d(f-k)\|_{L^2},
\end{equation}
where these $L^2$ norms use $\widehat g$ and $\mu(t_0)$. In particular, $E_t$ is continuous on $H^1(M)$. These arguments do not require reversibility.

\smallskip
\noindent\emph{Time dependence of the energy.}
Integrating the first equation in \eqref{F2} gives
\[
    e^{-A|t-s|}F_s(x,v)\leq F_t(x,v)
       \leq e^{A|t-s|}F_s(x,v).
\]
If $\mathbb B_t(x)=\{v\in T_xM\mid F_t(x,v)\leq1\}$, then
\[
    e^{-A|t-s|}\mathbb B_s(x)\subset\mathbb B_t(x)
       \subset e^{A|t-s|}\mathbb B_s(x).
\]
Taking support functions and using
$F_t^*(x,\xi)=\sup_{v\in\mathbb B_t(x)}\xi(v)$ therefore yields
\[
    e^{-A|t-s|}F_s^*(x,\xi)\leq F_t^*(x,\xi)
       \leq e^{A|t-s|}F_s^*(x,\xi).
\]
With $K=2A+B$, the metric and measure comparisons imply
\begin{equation}\label{eq:app-energy-time}
    e^{-K|t-s|}E_s(f)\leq E_t(f)
       \leq e^{K|t-s|}E_s(f)
       \qquad(f\in H^1(M)).
\end{equation}
Consequently,
\[
    |E_t(f)-E_s(f)|
       \leq\bigl(e^{K|t-s|}-1\bigr)E_t(f)
       \leq Ke^{K(b-t_0)}|t-s|E_t(f).
\]
This is \eqref{eq:app-relative-time}, with a constant independent of $f$. The function is held fixed in this estimate, so it is a property of the prescribed family of energies.

\smallskip
\noindent\emph{Identification of the subdifferential.}
For $u,\varphi\in H^1(M)$ and $0<|\varepsilon|\leq1$, the integral mean value formula gives the pointwise bound
\[
    \left|\frac{e_t(du+\varepsilon\,d\varphi)-e_t(du)}
                     {\varepsilon}\right|
       \leq C\bigl(|du|+|d\varphi|\bigr)|d\varphi|.
\]
The right-hand side is integrable. Differentiating under the integral, including on the set where $du=0$, gives
\begin{equation}\label{eq:app-first-variation}
    DE_t(u)[\varphi]
       =\left.\frac{d}{d\varepsilon}\right|_{\varepsilon=0}
                       E_t(u+\varepsilon\varphi)
       =\int_Md\varphi(\nabla u)\,d\mu(t).
\end{equation}
This is a directional derivative on $H^1(M)$.

If $\zeta\in\partial^{\,t}E_t(u)$, substitute $u+\varepsilon\varphi$ in \eqref{eq:app-subdifferential}, use both signs of $\varepsilon$, and let $\varepsilon\to0$. It follows that
\[
    \int_M\zeta\varphi\,d\mu(t)
       =\int_Md\varphi(\nabla u)\,d\mu(t)
       \qquad(\varphi\in H^1(M)).
\]
Thus $\zeta=-\Delta_{\mu(t)}u$ in the distributional sense. Conversely, if $\Delta_{\mu(t)}u\in L^2(M,\mu(t))$, convexity and \eqref{eq:app-first-variation} give
\[
    E_t(f)-E_t(u)\geq DE_t(u)[f-u]
       =\langle-\Delta_{\mu(t)}u,f-u\rangle_t
       \qquad(f\in H^1(M)).
\]
For $f\notin H^1(M)$ the subgradient inequality is automatic. Hence
\begin{equation}\label{eq:app-subgradient-laplacian}
\begin{aligned}
    \operatorname{Dom}(\partial^{\,t}E_t)
       &=\{u\in H^1(M)\mid\Delta_{\mu(t)}u\in L^2(M,\mu(t))\},\\
    \partial^{\,t}E_t(u)&=\{-\Delta_{\mu(t)}u\}.
\end{aligned}
\end{equation}
All the hypotheses of Theorem~\ref{thm:app-kopfer} have now been verified, and its subgradient evolution is exactly \eqref{eq:coupled-heat-cauchy}.

\smallskip
\noindent\emph{Implicit time approximation.}
We describe the construction and the bounds that give the asserted regularity. Choose
\[
    \delta=\frac{b-t_0}{N},\qquad
    t_j=t_0+j\delta,\qquad u^0=u_0,
\]
and successively minimize
\begin{equation}\label{eq:app-implicit-step}
    u^j=\operatorname*{argmin}_{f\in H^1(M)}
       \left\{E_{t_j}(f)
              +\frac1{2\delta}\|f-u^{j-1}\|_{t_j}^2\right\}.
\end{equation}
This is the time-dependent minimizing scheme in \cite[Section~5.1]{Kopfer2018GradientFlow}. The energy bounds $df$, and the quadratic term bounds $f$, so a minimizing sequence is bounded in $H^1(M)$. Weak compactness in $H^1(M)$, compactness of its inclusion into $L^2(M)$, and lower semicontinuity give a minimizer. The quadratic term makes the objective strictly convex, so the minimizer is unique. For the compact inclusion one may use the auxiliary Riemannian metric, as in \cite[Lemma~2.2(ii)]{OhtaSturm2009HeatFlow}.

The Euler equation for \eqref{eq:app-implicit-step} is
\begin{equation}\label{eq:app-discrete-heat}
    \frac1\delta\int_M(u^j-u^{j-1})\varphi\,d\mu(t_j)
       +\int_Md\varphi\bigl(\nabla^{F_{t_j}}u^j\bigr)\,d\mu(t_j)=0,
       \qquad\varphi\in H^1(M).
\end{equation}
In particular, each step solves the weak implicit heat equation at the current metric and measure. Comparing the minimizer with $u^{j-1}$ and using \eqref{eq:app-energy-time}, we obtain
\[
    E_{t_j}(u^j)
       +\frac{\|u^j-u^{j-1}\|_{t_j}^2}{2\delta}
       \leq e^{K\delta}E_{t_{j-1}}(u^{j-1}).
\]
Multiplication by $e^{-K(t_j-t_0)}$ and summation give
\begin{equation}\label{eq:app-discrete-energy}
\begin{aligned}
    &e^{-K(t_m-t_0)}E_{t_m}(u^m)\\
    &\quad+\sum_{j=1}^m e^{-K(t_j-t_0)}
              \frac{\|u^j-u^{j-1}\|_{t_j}^2}{2\delta}
       \leq E_{t_0}(u_0),\qquad 1\leq m\leq N.
\end{aligned}
\end{equation}
It follows that
\[
\begin{aligned}
    \max_{0\leq j\leq N}E_{t_j}(u^j)
       &\leq e^{K(b-t_0)}E_{t_0}(u_0),\\
    \sum_{j=1}^N\delta
       \left\|\frac{u^j-u^{j-1}}{\delta}\right\|_{t_j}^2
       &\leq 2e^{K(b-t_0)}E_{t_0}(u_0).
\end{aligned}
\]
These are the energy and velocity bounds underlying \cite[Proposition~5.3]{Kopfer2018GradientFlow}.

\smallskip
\noindent\emph{Passage to the limit.}
For $t\in(t_{j-1},t_j]$, define
\[
    \widehat u_\delta(t)
       =u^{j-1}+\frac{t-t_{j-1}}{\delta}(u^j-u^{j-1}),
       \qquad
    \overline u_\delta(t)=u^j,
\]
and set both interpolants equal to $u_0$ at $t_0$. The preceding bounds and \eqref{eq:app-l2-time} imply a uniform bound for $\partial_t\widehat u_\delta$ in $L^2((t_0,b);\mathcal H)$. Integration from the fixed initial value bounds $\widehat u_\delta$ in $L^\infty((t_0,b);\mathcal H)$. The spatial energy bound then gives a uniform $L^\infty((t_0,b);H^1(M))$ bound for both interpolants. Moreover,
\[
\begin{aligned}
    \|\widehat u_\delta(t)-\widehat u_\delta(s)\|_{\mathcal H}
       &\leq C|t-s|^{1/2},\\
    \sup_{t\in[t_0,b]}
       \|\overline u_\delta(t)-\widehat u_\delta(t)\|_{\mathcal H}
       &\leq C\delta^{1/2}.
\end{aligned}
\]
The compact inclusion $H^1(M)\hookrightarrow L^2(M)$ and the Arzel\`a--Ascoli theorem therefore give a subsequence such that
\begin{equation}\label{eq:app-interpolant-convergence}
\begin{aligned}
    \widehat u_\delta&\longrightarrow u
        &&\text{in }C([t_0,b];\mathcal H),\\
    \overline u_\delta&\longrightarrow u
        &&\text{uniformly on }[t_0,b]\text{ in }\mathcal H,\\
    \partial_t\widehat u_\delta&\rightharpoonup u_t
        &&\text{in }L^2((t_0,b);\mathcal H).
\end{aligned}
\end{equation}
This is the convergence supplied in the general Hilbert-space setting by \cite[Proposition~5.4]{Kopfer2018GradientFlow}. Here compactness of $M$ also gives the direct argument above. Lower semicontinuity, the time comparison of the energies, and \eqref{eq:app-discrete-energy} imply
\[
    u\in H^1((t_0,b);\mathcal H)
          \cap L^\infty((t_0,b);H^1(M)),
    \qquad u(t_0)=u_0.
\]

To identify the equation, use the variational inequality for each minimizer. For $f\in H^1(M)$ and $t\in(t_{j-1},t_j)$, it reads
\[
    \langle\partial_t\widehat u_\delta,
                \overline u_\delta-f\rangle_{t_j}
       +E_{t_j}(\overline u_\delta)-E_{t_j}(f)\leq0.
\]
The smooth positive densities satisfy
\[
    |\langle v,w\rangle_t-\langle v,w\rangle_s|
       \leq C|t-s|\|v\|_{\mathcal H}\|w\|_{\mathcal H}.
\]
Together with \eqref{eq:app-relative-time} and the uniform energy and velocity bounds, this permits replacement of $t_j$ by $t$ with an integrated error tending to zero. Multiply by a nonnegative $\eta\in C_c^\infty((t_0,b))$ and integrate. The scalar-product term converges by the weak convergence of the velocities and the strong convergence of $\overline u_\delta$. Lower semicontinuity and Fatou's lemma apply to the nonnegative energy term. We obtain
\begin{equation}\label{eq:app-limit-inequality}
    \int_{t_0}^b\eta(t)
       \left[\langle u_t,u-f\rangle_t
               +E_t(u)-E_t(f)\right]\,dt\leq0.
\end{equation}
First take $f$ in a countable dense subset of $H^1(M)$. The fundamental lemma for nonnegative time tests gives one exceptional null set outside which the integrand is nonpositive for every such $f$. By \eqref{eq:app-h1-continuity} the inequality extends to every $f\in H^1(M)$. For $f\notin H^1(M)$ it is automatic. Thus
\[
    -u_t\in\partial^{\,t}E_t(u)
    \quad\text{for almost every }t.
\]
The identification \eqref{eq:app-subgradient-laplacian} gives \eqref{eq:coupled-heat-cauchy} and its stated weak formulation. In particular, the distributional Laplacian is represented by the $L^2$ function $u_t$ for almost every time, although this is not required of $u_0$.

\smallskip
\noindent\emph{Uniqueness and nonnegativity.}
If $u$ and $\widetilde u$ are two solutions in the asserted class and $w=u-\widetilde u$, the chain rule in $L^2$ and the measure evolution give
\[
\begin{aligned}
    \frac12\frac{d}{dt}\int_Mw^2\,d\mu(t)
       ={}&-\int_Mdw(\nabla u-\nabla\widetilde u)\,d\mu(t)
             -\frac12\int_M\mathcal Rw^2\,d\mu(t)\\
       \leq{}&\frac B2\int_Mw^2\,d\mu(t).
\end{aligned}
\]
Indeed, convexity of $e_t$ implies
\[
    (\xi-\zeta)\bigl(J_t^*(\xi)-J_t^*(\zeta)\bigr)\geq0,
\]
which gives the sign of the first integral. Gronwall's inequality yields
\begin{equation}\label{eq:app-uniqueness}
    \|u(t)-\widetilde u(t)\|_t^2
       \leq e^{B(t-t_0)}
                \|u(t_0)-\widetilde u(t_0)\|_{t_0}^2.
\end{equation}
Hence the solution is unique.

Nonnegativity follows already at the discrete level. If $u^{j-1}\geq0$, replace a competitor $f$ in \eqref{eq:app-implicit-step} by $f_+=\max\{f,0\}$. The Sobolev chain rule gives $df_+=\mathbf1_{\{f>0\}}df$, and therefore
\[
    E_{t_j}(f_+)\leq E_{t_j}(f),\qquad
    \|f_+-u^{j-1}\|_{t_j}\leq\|f-u^{j-1}\|_{t_j}.
\]
Uniqueness of the minimizer implies $u^j\geq0$. Induction and the strong $L^2$ limit show that $u\geq0$ whenever $u_0\geq0$. This argument also applies when $F_t$ is nonreversible.

Finally, testing the weak equation with $u$ and differentiating the measure gives
\begin{equation}\label{eq:app-l2-energy-identity}
    \frac12\frac{d}{dt}\int_Mu^2\,d\mu(t)
       +\int_MF_t^{*2}(du)\,d\mu(t)
       =-\frac12\int_M\mathcal R u^2\,d\mu(t)
\end{equation}
for almost every time. The term involving $\mathcal R$ comes from the derivative of the measure. The subgradient is taken at fixed time, so the scalar evolution itself is \eqref{eq:coupled-heat-cauchy}. This completes the proof.
\end{proof}

The argument applies on every compact interval on which the prescribed metric and measure remain smooth and the metric remains strongly convex. It gives the variational regularity in \eqref{eq:coupled-heat-existence-regularity}. The additional spatial regularity required by the frequency identities remains the assumption in \eqref{eq:independent-regularity}. Local smooth solutions on a noncritical region are supplied separately by Proposition~\ref{prop:local-smooth-heat}.

\section*{Acknowledgments}

The authors are grateful to the reviewers for their careful reading and valuable comments.

\end{document}